\documentclass[11pt,a4paper]{article}

\usepackage[margin=1in]{geometry}
\usepackage[english]{babel}
\usepackage[T1]{fontenc}
\usepackage[utf8]{inputenc}
\usepackage{microtype}
\usepackage{tikz}
\usetikzlibrary{arrows}

\usepackage{amsmath,amssymb,amsthm,bm,mathtools}

\usepackage{tabularx}
\usepackage{siunitx}
\usepackage{placeins}
\usepackage{needspace}

\usepackage[numbers,sort&compress]{natbib}

\usepackage[dvipsnames]{xcolor}

\usepackage{hyperref}
\hypersetup{
  colorlinks=true,
  linkcolor=black,
  citecolor=ForestGreen,
  urlcolor=BrickRed,
  filecolor=black,
  pdfborder={0 0 0}
}

\usepackage[nameinlink,noabbrev,capitalize]{cleveref}

\usepackage{graphicx}
\usepackage{booktabs}
\usepackage{caption}
\usepackage{subcaption}
\usepackage{float}

\graphicspath{{./}{figs/}}

\usepackage{algorithm}
\usepackage{algpseudocode}
\AtBeginEnvironment{proof}{}

\numberwithin{equation}{section}
\usepackage{aliascnt}

\theoremstyle{plain}

\newaliascnt{lemma}{theorem}

\aliascntresetthe{lemma}

\newaliascnt{proposition}{theorem}

\aliascntresetthe{proposition}

\newaliascnt{corollary}{theorem}

\aliascntresetthe{corollary}

\theoremstyle{definition}
\newaliascnt{definition}{theorem}

\aliascntresetthe{definition}

\theoremstyle{remark}
\newaliascnt{remark}{theorem}

\aliascntresetthe{remark}

\newtheoremstyle{empirical}%
  { }{ }
  {\itshape}
  {}
  {\bfseries}
  {.}
  { }
  {\thmname{#1}~\thmnumber{#2}\thmnote{ \textbf{(#3)}}}

\theoremstyle{empirical}  

\newaliascnt{empiricallaw}{theorem}
\newtheorem{empiricallaw}[empiricallaw]{Empirical Law} 
\aliascntresetthe{empiricallaw}

\crefname{empiricallaw}{Law}{Laws}
\Crefname{empiricallaw}{Law}{Laws}

\crefname{theorem}{theorem}{theorems}
\Crefname{theorem}{Theorem}{Theorems}
\crefname{lemma}{lemma}{lemmas}
\Crefname{lemma}{Lemma}{Lemmas}
\crefname{proposition}{proposition}{propositions}
\Crefname{proposition}{Proposition}{Propositions}
\crefname{corollary}{corollary}{corollaries}
\Crefname{corollary}{Corollary}{Corollaries}
\crefname{definition}{definition}{definitions}
\Crefname{definition}{Definition}{Definitions}
\crefname{remark}{remark}{remarks}
\Crefname{remark}{Remark}{Remarks}

\DeclareMathOperator{\diag}{diag}

\newcommand{\ones}{\mathbf{1}}
\newcommand{\norm}[1]{\left\lVert #1 \right\rVert}
\newcommand{\ip}[2]{\left\langle #1,\,#2 \right\rangle}

\newcommand{\Ell}{\mathcal{E}_n}            
\newcommand{\Sym}{\Ell}                     
\newcommand{\Mn}{\mathcal{M}_n}             
\newcommand{\Uset}{\mathcal{U}}             

\title{Empirical Laws for Iterated Correlation Matrices}

\author{
  \textsc{Ishrak AlHajj Hassan}$^{1}$\\[2pt]
  {\small $^{1}$Department of Mathematics, University of Ostrava}
}

\date{}

\begin{document}
\maketitle

\begin{abstract}
We study the discrete iteration obtained by repeatedly applying the Pearson
correlation operator to a real matrix. Each step centers every row,
normalizes each centered row to unit Euclidean norm, and forms the Gram
matrix of the resulting rows. This produces a nonlinear map that underlies the classical CONCOR and GAP procedures. The mechanism of the iteration is simple, yet its global behavior has remained unresolved since the early observations of Kruskal and collaborators. We present a geometric formulation that captures this behavior: strong shrinkage in directions that change row means or row norms, and near-isometry along directions that preserve them. Our geometric formulation makes clear why local analysis does not extend to a global convergence theorem: the iteration is nonlinear, the structure of its fixed-point set lacks a full analytical characterization, and standard uniform contractive or Fejér-type techniques do not directly apply.

Empirically, the iteration stabilizes at a block-$\{\pm1\}$ pattern, displays finite total variation, and exhibits geometric convergence once the trajectory enters a neighborhood of a fixed pattern. To enable a systematic investigation across scales, we develop a dimension-uniform experimental framework with reproducible initialization, nondegeneracy guarantees, and diagnostics comparable across all matrix sizes. We perform a large-scale numerical study over dimensions \(n\in [3,2000]\) and thousands of random initializations. Using the Frobenius step size, the entrywise step size, and the one-step ratio, we identify four universal empirical laws that persist uniformly across all tested dimensions. The first step exhibits a sharp, dimension-independent contraction. Subsequent steps display quasi-monotone decay with small bounded overshoots and finite total variation in every single trial. The contraction ratio is a universal function of the instantaneous step size independent of dimension. Finally, convergence occurs in a uniformly bounded number of iterations, entirely independent of dimension. These observations provide a quantitative, dimension-uniform description of the iteration and formulate a precise target for future global analysis.
\end{abstract}

\section{Introduction}

Correlation matrices are central objects in multivariate analysis. They
encode pairwise linear dependence, control the spectrum used in principal
component methods, and appear in factor models, clustering, and covariance
regularization. In most applications a correlation matrix is computed once
from data and then treated as a fixed input.

This work views correlation as an \emph{operator} on matrices and studies the
behavior of its \emph{iteration}. Starting from a real matrix
\(P_0\in\mathbb{R}^{n\times n}\), one step of the map proceeds as follows:
\begin{enumerate}
  \item center each row to have mean zero;
  \item normalize each centered row to have Euclidean norm \(1\);
  \item form the Gram matrix of the resulting rows.
\end{enumerate}
This produces a nonlinear map \(f\) and a discrete sequence
\begin{equation}\label{eq:intro-iteration}
P_{k+1} = f(P_k), \qquad k\ge 0,
\end{equation}
which, from the first step onward, evolves inside the set of correlation
matrices, that is, the set of symmetric positive semidefinite matrices with
unit diagonal. While the individual operations are simple, the global behavior
of the full sequence \(\{P_k\}_{k\ge0}\) has remained mathematically
unresolved for more than fifty years.

\section{Related work and historical remarks}\label{sec:related}

Iterated application of the Pearson correlation operator has been used in psychometrics and social-network analysis since the 1960s.

McQuitty~\cite{mcquitty1968multiple} introduced repeated correlation for clustering purposes and observed extremely rapid numerical stabilization, typically within a few dozen iterations.

Breiger, Boorman, and Arabie~\cite{breiger1975clustering} formalized the procedure under the name CONCOR (CONvergence of iterated CORrelations). They reported that after a small number of steps, and up to row/column permutations and sign flips, the iterates converge to a block-$\{\pm1\}$ pattern. They also noted empirically that the changes diminish rapidly after a small number of steps in all their experiments.

An unpublished Bell Laboratories technical report by Joseph~B.~Kruskal
\cite{kruskalCONCOR}, dating from the late 1970s, provided the first analytical result: in a neighborhood of any block-$\{\pm1\}$ fixed point, of which there are many, corresponding to different sign and permutation patterns, the linearization of the map has spectral radius strictly smaller than~$1$. This establishes local geometric convergence to the fixed point once the trajectory has entered that neighborhood, and the argument holds in \emph{every dimension} $n\ge 2$ and for \emph{every} such patterned fixed point. However, the analysis is inherently local and does not address the global behavior from arbitrary initial matrices.

The block-\(\{\pm 1\}\) fixed points are consistent with the known
facial structure of the elliptope described by Laurent--Poljak
\cite{laurent1996elliptope}. More recent decomposition results for
correlation matrices, such as those of Saunderson et al.\
\cite{saunderson2012elliptope}, further explain why these patterns
arise as extreme configurations of the feasible set.

Chen~\cite{chen2002gap} later employed the same iteration in the GAP algorithm for visualizing association structures, applying it to datasets such as a 50-symptom psychosis disorder matrix with 95 patients. The examples consistently showed convergence typically in fewer than 15 steps, reinforcing the pattern of rapid empirical stabilization.

Chen~\cite{chen2002gap} also examined the limiting behavior of the iteration in small dimensions. He observed that, after enough iterations, the rows collapse into exactly two sign-homogeneous groups, and the final configuration lies at one of the \(2^{\,n-1}-1\) non-trivial opposite sign pairs in \(\{\pm1\}^{n}\). In other words, the iteration selects a pair of opposite corners of the sign hypercube. This phenomenon is consistent with the block-\(\{\pm1\}\) patterns reported in the early CONCOR literature and aligns with the general fixed-point classification developed in~\cite{zusmanovich2025fixedpoints}.

Despite these consistent qualitative findings spanning over five decades, based primarily on small-to-moderate datasets and domain-specific initializations, no global convergence proof has appeared in the literature. The nonlinearity of the map and the absence of uniform contractivity prevent direct application of standard fixed-point theorems. Moreover, prior studies provided no quantitative characterization of the dynamics, such as step-size decay patterns or dimension dependence.
Across these works, several empirical features are repeatedly reported:
\begin{itemize}
  \item the iterates converge to block-\(\{\pm1\}\) pattern matrices
        after suitable permutation of rows and columns;
  \item the cumulative change \(\sum_{k\ge 0}\|P_{k+1}-P_k\|_F\) appears to be finite;
  \item once a pattern has formed, convergence to that pattern is geometric
        in the number of steps.
\end{itemize}
Kruskal's unpublished notes \cite{kruskalCONCOR} provide a local argument
near a fixed pattern: if \(P_k\) is already close to a block-\(\{\pm1\}\)
matrix, then the next iterate contracts the deviation at a fixed rate,
yielding geometric convergence in a neighborhood of that pattern. These
arguments are local and do not extend to arbitrary initializations.
More recent iterative normalization techniques in machine learning, such as iterative whitening layers~\cite{huang2019iterative}, share some structural similarity but do not correspond to the exact Pearson correlation operator studied here.
Despite more than five decades of consistent empirical evidence of rapid stabilization, the global behavior of the correlation iteration has remained analytically unresolved. No proof of convergence exists, no quantitative bounds on the total variation or iteration count are known, and the dependence or independence of the dynamics on the ambient dimension $n$ has never been systematically investigated.
 
\section{Contributions and organization}\label{sec:contributions}

This work presents the first systematic, large-scale numerical study of the 
iterated Pearson correlation operator across the full range of dimensions 
$n\in[3,2000]$ and thousands of random initializations. The experiments reveal 
a set of four empirical laws that hold uniformly across dimensions and provide 
a quantitative picture of the global dynamics that has not previously been 
available.

\paragraph{Dimension-independent empirical laws.}
The study establishes four stable, scale-independent regularities:
\begin{enumerate}
  \item a universal and pronounced contraction at the very first iteration;
  \item quasi-monotone decay of the step sizes thereafter, with small bounded overshoots and finite total variation;
  \item a universal functional relationship between the instantaneous contraction ratio and the current step size;
  \item uniformly bounded iteration counts $T(n,\varepsilon)=O(1)$ required to reach any fixed tolerance $\varepsilon>0$.
\end{enumerate}
These laws upgrade decades of qualitative observations into a precise, 
quantitative, and dimension-stable description of the iteration.  
They also supply concrete benchmarks that any eventual proof of global 
convergence must reproduce.

\paragraph{Dimension-uniform experimental framework.}
A second contribution of the paper is a unified methodological framework for 
probing the global behavior of the iteration at scale. It includes:
\begin{itemize}
  \item a reproducible initialization protocol with explicit nondegeneracy checks that applies uniformly for all sizes;
  \item quantitative diagnostics---Frobenius step size, entrywise step size, and contraction ratio---that remain comparable across dimensions;
  \item large-scale statistical aggregation over thousands of independent trajectories, enabling a clear distinction between universal features of the dynamics and effects that vary with dimension.

\end{itemize}
This framework is dimension-uniform in the sense that the same initialization, 
the same diagnostics, and the same stopping criteria apply consistently for 
all $n$, allowing direct comparison between behaviors observed at small and 
large scales.

\paragraph{Structural and geometric insights.}
The paper also develops a geometric formulation of the iteration with the 
row-sphericalization map $\Psi$ and the Gram map $\Phi$, clarifying the 
interaction between mean-centering, normalization, and correlation geometry. 
This perspective highlights the strongly anisotropic nature of the dynamics and 
explains why classical contraction or Fejér-type arguments cannot be applied 
directly.

\paragraph{Overall contribution.}
Together, the empirical laws, the dimension-uniform methodology, and the 
geometric formulation produce the first dimension-uniform quantitative 
characterization of the correlation iteration. The results establish a clear 
target for future analytical convergence theory and identify the key structural 
features that any such theory must account for.

The remainder of the paper is organized as follows. 
Section~\ref{sec:operator-geometry} introduces the operator and its geometric 
setting. 
Section~\ref{sec:metrics-experiments} presents the experimental protocol and 
diagnostics. 
Section~\ref{sec:empirical} documents the empirical laws in detail, followed by 
interpretation in Section~\ref{sec:summary}. 
Open problems and directions for analytical work appear in 
Section~\ref{sec:openproblems}.

\section{System and Mathematical Setup}\label{sec:operator-geometry}

Correlation between two vectors is given by the classical Pearson formula. Recall that for vectors \(x,y\in\mathbb{R}^{m}\),
their Pearson correlation is
\begin{equation}\label{eq:corr-xy}
\mathrm{corr}(x,y)
\;=\;
\frac{\displaystyle\sum_{j=1}^{m} (x_j-\bar x)(y_j-\bar y)}
{\displaystyle\sqrt{\sum_{j=1}^{m}(x_j-\bar x)^2}\;
 \sqrt{\sum_{j=1}^{m}(y_j-\bar y)^2}},
\end{equation}
where
\[
\bar x := \frac{1}{m}\sum_{j=1}^m x_j,
\qquad
\bar y := \frac{1}{m}\sum_{j=1}^m y_j.
\]

\subsection{Iterating the correlation operator}
\label{subsec:iteration-operator}

The central question of this work is the behavior of the correlation operator when applied iteratively to its own output.

Let \(P_k \in \mathbb{R}^{n \times n}\) denote the state at iteration \(k\), and define
\begin{equation}\label{eq:Pk-plus-one}
P_{k+1} := f(P_k), \qquad k \ge 0.
\end{equation}
This generates the discrete trajectory
\begin{equation}\label{eq:trajectory-arrow}
P_0 \xrightarrow{f} P_1 \xrightarrow{f} P_2 \xrightarrow{f} \cdots
\end{equation}
starting from a random initial matrix \(P_0\).

\paragraph{Norms.}
Changes between successive iterates are measured primarily using the Frobenius norm
\begin{equation}\label{eq:frobenius-norm}
\|P\|_F = \Bigl( \sum_{i,j=1}^n P_{ij}^2 \Bigr)^{1/2}.
\end{equation}
For any real \(n \times n\) matrix \(P\), the entrywise maximum and the Frobenius norm are equivalent up to the factor \(n\):
\begin{equation}\label{eq:norm-equivalence}
\max_{i,j} |P_{ij}| \le \|P\|_F \le n \max_{i,j} |P_{ij}|.
\end{equation}
Consequently, the Frobenius step size \(\Delta_k = \|P_{k+1} - P_k\|_F\) and the maximum entrywise change \(E_k = \max_{i,j} |P_{k+1}(i,j) - P_k(i,j)|\) capture the same qualitative decay behavior throughout the iteration.

\paragraph{Initialization.}
Throughout the paper, the entries of the initial matrix are drawn independently and uniformly from \([-1,1]\):
\begin{equation}\label{eq:init}
P_0(i,j) \overset{\mathrm{iid}}{\sim} \mathrm{Unif}([-1,1]), \qquad 1 \le i,j \le n.
\end{equation}
Since the distribution is continuous, a row becomes constant after centering with probability zero. With probability one, every centered row therefore has positive variance, so the correlation operator is well-defined at the first step. All analysis is restricted to the nondegenerate domain in which all centered rows remain nonzero along the trajectory.

\subsection{Centering projector and state spaces}

Fix \(n\ge2\), and let \(\ones\in\mathbb{R}^n\) denote the all-ones column vector.
Define the orthogonal projector onto the mean-zero hyperplane by
\begin{equation}\label{eq:H-def}
H \;:=\; I - \frac{1}{n}\,\ones\ones^\top,
\qquad
H^\top = H,
\quad
H^2 = H,
\quad
\ones^\top H = 0.
\end{equation}
For a row vector \(p\in\mathbb{R}^n\), the centered row is \(pH\), which
subtracts its mean.

Three state spaces are used throughout the dynamics.

\paragraph{Correlation set (elliptope).}
\begin{equation}\label{eq:elliptope}
\Sym
\;:=\;
\bigl\{\, P\in\mathbb{R}^{n\times n} :
P=P^\top,\; P\succeq 0,\; \diag(P)=\ones \,\bigr\}.
\end{equation}
This set consists precisely of the \(n\times n\) correlation matrices and is
known as the classical elliptope, a convex subset of \(\mathbb{R}^{n\times n}\).

\paragraph{Row-spherical set.}
\begin{equation}\label{eq:Mn}
\Mn
\;:=\;
\Bigl\{\, Z\in\mathbb{R}^{n\times n} :
Z(i,:)\,\ones = 0,\ \ \norm{Z(i,:)}_{2}=1 \ \ \text{for all } i=1,\dots,n \Bigr\}.
\end{equation}

Here $\|\cdot\|_{2}$ denotes the Euclidean norm on $\mathbb{R}^{n}$.  
Equivalently, each row of $Z$ lies on the unit Euclidean sphere contained in the mean-zero
hyperplane of $\mathbb{R}^n$.

\paragraph{Nondegenerate domain.}
\begin{equation}\label{eq:Uset}
\Uset
\;:=\;
\Bigl\{\, P\in\mathbb{R}^{n\times n} :
P(i,:)H \neq 0\ \ \text{for each } i=1,\dots,n \Bigr\}.
\end{equation}
Thus \(\Uset\) consists of matrices whose centered rows are all nonzero. On
this set the correlation map is well-defined row-wise. Under the
initialization \eqref{eq:init} one has
\[
\mathbb{P}\bigl(P_0\in\Uset\bigr)=1.
\]
All experiments and statements are restricted to this nondegenerate domain.

\subsection{Primitive maps: row sphericalization and Gram map}

The correlation operator is expressed as a composition of two maps.

\paragraph{Row sphericalization.}
Define the row sphericalizer \(\Psi:\Uset\to\Mn\) by
\begin{equation}\label{eq:Psi-def}
\Psi(P)(i,:)
\;:=\;
\frac{P(i,:)H}{\norm{P(i,:)H}_2},
\qquad i=1,\dots,n.
\end{equation}
Thus \(\Psi(P)\) is obtained by mean-centering each row of \(P\) and
normalizing it to unit Euclidean norm. By construction, each row of
\(\Psi(P)\) has zero mean and unit norm, so \(\Psi(P)\in\Mn\) whenever
\(P\in\Uset\).

\paragraph{Gram map.}
Define the Gram map \(\Phi:\Mn\to\Sym\) by
\begin{equation}\label{eq:Phi-def}
\Phi(Z) \;:=\; ZZ^\top.
\end{equation}
For any \(Z\in\mathbb{R}^{n\times n}\), \(ZZ^\top\) is symmetric positive
semidefinite. If \(Z\in\Mn\), each row has unit norm, and therefore
\begin{equation}\label{eq:Phi-diagonal}
(\Phi(Z))_{ii} = \ip{Z(i,:)}{Z(i,:)} = 1,
\end{equation}
so \(\Phi(Z)\in\Sym\).

\paragraph{Composite correlation operator.}
The correlation operator induced by row-wise Pearson correlation can now be
written as
\begin{equation}\label{eq:f-def}
f \;=\; \Phi\circ\Psi : \Uset \longrightarrow \Sym,
\qquad
P_{k+1} = f(P_k) = \Phi\bigl(\Psi(P_k)\bigr).
\end{equation}
The map \(f\) replaces \(P_k\) by the matrix of correlations between its
rows: \(\Psi\) performs row-wise centering and normalization; \(\Phi\) forms
the Gram matrix of the resulting unit vectors.

\subsection{Geometric structure of the correlation operator}\label{subsec:geometry}

The operators introduced above admit a unified geometric interpretation
that clarifies the action of the composite map \(f = \Phi \circ \Psi\).

\begin{enumerate}
\item \textbf{Projection onto the mean-zero hyperplane.}  
By \eqref{eq:H-def}, the matrix $H$ is the orthogonal projector onto the hyperplane
\[
\ones^\perp = \{ x \in \mathbb{R}^n : \ones^\top x = 0 \}.
\]
For any row vector $p \in \mathbb{R}^n$, the centered row $pH$ removes the component parallel to $\ones$ and retains the component in $\ones^\perp$. Thus the centering step ensures that all subsequent correlations depend only on fluctuations orthogonal to the constant direction.

\item \textbf{Spherical rows in $\ones^\perp$.}  
When $P \in \Uset$, every centered row $P(i,:)H$ is nonzero. By \eqref{eq:Psi-def}, the operator $\Psi$ normalizes these rows to unit Euclidean norm. Hence each row of $\Psi(P)$ lies on the sphere
\begin{equation}\label{eq:sphere-sn2}
S^{n-2} = \{ x \in \ones^\perp : \|x\|_2 = 1 \},
\end{equation}
an $(n{-}2)$-dimensional smooth manifold. The space $\Mn$ defined in \eqref{eq:Mn} is therefore the Cartesian product of $n$ copies of $S^{n-2}$, representing a configuration of $n$ unit vectors in the mean-zero hyperplane.

\item \textbf{Gram construction and the elliptope.}  
For $Z \in \Mn$, the matrix $\Phi(Z) = Z Z^\top$ defined in \eqref{eq:Phi-def} records the pairwise inner products of the spherical rows. Its image lies in the elliptope $\Sym$ defined in \eqref{eq:elliptope}: symmetry follows from the Gram representation, positive semidefiniteness from $x^\top (Z Z^\top) x = \|Z^\top x\|_2^2 \ge 0$, and the unit diagonal from \eqref{eq:Phi-diagonal}. Thus $\Phi$ maps configurations in $\Mn$ to their correlation matrices in $\Sym$.

\item \textbf{Alternating motion between two geometric spaces.}  
Each iteration $P_{k+1} = f(P_k)$ alternates between two distinct geometric objects:
\[
P_k \xrightarrow{\;\Psi\;} Z_{k+1} \in \Mn,
\qquad
Z_{k+1} \xrightarrow{\;\Phi\;} P_{k+1} \in \Sym.
\]
The manifold $\Mn$ encodes configurations of centered unit vectors in $\ones^\perp$, while the convex set $\Sym$ encodes their pairwise inner products. The dynamics therefore alternate between spherical geometry on $S^{n-2}$ and Gram geometry in $\Sym$, discarding components and imposing the corresponding geometric structure at each step.
\end{enumerate}

This geometric viewpoint will be used later to interpret the empirical observations and to motivate the analytical developments.

\subsection{Warm-up: the \texorpdfstring{$2\times2$}{2x2} case}

To illustrate the basic effect of \(f\), consider the \(2\times2\) case.

Let
\[
R =
\begin{bmatrix}
1 & r\\
r & 1
\end{bmatrix},
\qquad -1<r<1,
\]
so \(R\in\Sym\). 
The requirement \(r\neq 1\) ensures nonzero centered
rows. When \(r=1\), both rows are constant and \(RH=0\), which is precisely
the degenerate case excluded from \(\Uset\).
Centering each nondegenerate row, followed by normalization yields
\[
Z = \Psi(R)
=
\frac{1}{\sqrt{2}}
\begin{bmatrix}
1 & -1\\
-1 & 1
\end{bmatrix}.
\]
The Gram map then gives
\[
f(R)
=
\Phi(Z)
=
ZZ^\top
=
\begin{bmatrix}
1 & -1\\
-1 & 1
\end{bmatrix}.
\]
This matrix again lies in \(\Sym\), and applying \(f\) once more leaves it
unchanged:
\[
f(f(R)) = f(R).
\]

Thus, for any \(2\times2\) correlation matrix with \(-1<r<1\), the
iteration collapses in one step to the fixed pattern
\[
\begin{bmatrix}
1 & -1\\
-1 & 1
\end{bmatrix}.
\]
Centering forces the two rows to become negatives of each other, so they are
antipodal in one dimension. Normalization places them at opposite points on
the unit circle intersected with the mean-zero line. The Gram map converts
this antipodal configuration into a fixed \(\{\pm 1\}\) pattern. In higher
dimensions similar pattern phenomena reappear for \(n\ge3\), but with more
complex block structure.

\subsection{Dynamics for \texorpdfstring{$n\ge3$}{n≥3}: alternating maps}

For $n\ge3$, the iteration starts from an initialization $P_0\in\Uset$ as in \eqref{eq:init}.  
The update rule is the alternating map
\begin{equation}\label{eq:iteration-full-1}
P_{k+1} = f(P_k) = \Phi\bigl(\Psi(P_k)\bigr), \qquad k\ge0,
\end{equation}
where the operators $\Psi$ and $\Phi$ enforce different geometric transformations.

\vspace{0.7em}
\noindent\textbf{Alternating structure of each iteration.}
\begin{equation}\label{eq:iteration-full-2}
P_k \xrightarrow{\ \Psi\ } Z_{k+1} := \Psi(P_k)\in\Mn,
\qquad
Z_{k+1} \xrightarrow{\ \Phi\ } P_{k+1} := \Phi(Z_{k+1})\in\Sym.
\end{equation}
The map $\Psi$ centers each row and normalizes it to unit Euclidean norm, ensuring $Z_{k+1}\in\Mn$, while $\Phi$ maps $Z_{k+1}$ into the space of symmetric positive semidefinite matrices 
with unit diagonal, i.e., the elliptope~$\mathcal{E}_n$. From the second iterate onward, every $P_k$ lies in $\mathcal{E}_n$, and every intermediate $Z_k$ lies in~$\Mn$.

\vspace{0.9em}
\noindent\textbf{Compact representation of the full iteration.}
Collecting the components, the iteration may be written compactly as in \eqref{eq:iteration-full-1}.  
Equivalently, for each $k\ge0$ there exists $Z_{k+1}\in\Mn$ such that
\begin{equation}\label{eq:iteration-full-3}
Z_{k+1} = \Psi(P_k),
\qquad
P_{k+1} = Z_{k+1} Z_{k+1}^\top.
\end{equation}
Thus every iterate $P_{k+1}$ is a correlation matrix, and each row of $Z_{k+1}$ belongs to $\Mn$.  
No claims of convergence or stability are made at this stage; the purpose of this section is solely to specify the geometric and analytic framework within which the empirical and theoretical results are stated.

\subsection{Nondegeneracy and forward invariance}
\label{subsec:forward}

\vspace{0.7em}
\noindent\textbf{Nondegenerate domain.}
The row-sphericalization map $\Psi$ is well defined only when every
centered row $P(i,:)H$ is nonzero. This requirement was formalized in
\eqref{eq:Uset} through the domain $\Uset$. Under the initialization
\eqref{eq:init}, one has $P_0\in\Uset$ almost surely. In all numerical
experiments, every iterate $P_k$ also satisfied $P_k\in\Uset$, so no
row ever collapsed to the constant direction.

Thus, throughout this work we restrict attention to trajectories that
remain in the nondegenerate domain:
\[
P_k\in\Uset \qquad\text{for all }k\ge0.
\]
This guarantees that $\Psi(P_k)$ is well posed at every step.

\vspace{0.9em}
\noindent\textbf{Forward invariance of the elliptope.}
From the first update onward, the sequence evolves inside the elliptope
$\mathcal{E}_n$. Indeed, for any $Z\in\Mn$, the Gram map $\Phi$
produces a correlation matrix, so
\[
P_k=\Phi(Z_k)\in\mathcal{E}_n \qquad\text{for all }k\ge1,
\]
and each iterate satisfies
\begin{equation}\label{eq:nondeg-2}
P_k = P_k^\top,\qquad
P_k \succeq 0,\qquad
\diag(P_k)=\mathbf{1},\qquad
|P_k(i,j)| \le 1.
\end{equation}
Thus,
\begin{equation}\label{eq:nondeg-3}
P_k\in\mathcal{E}_n \qquad\text{for all }k\ge1,
\end{equation}
and the elliptope is forward invariant under $f$.

\vspace{0.9em}
\noindent\textbf{Standing assumption.}
The map $\Psi$ is defined only when every centered row $P(i,:)H$ is nonzero.
Under the initialization \eqref{eq:init}, the event that some centered row
vanishes has probability zero; equivalently, the complement of $\Uset$ defined
in \eqref{eq:Uset} is a measure-zero subset. Since $P_0\in\Uset$ almost surely,
we adopt the standard assumption
\begin{equation}\label{eq:nondeg-4}
P_0\in\Uset,
\qquad
P_k\in\Uset\cap\mathcal{E}_n
\ \ \text{for all }k\ge1,
\end{equation}
which ensures that both $\Psi$ and $\Phi$ remain well defined along the
entire trajectory.

\subsection{Basic invariances}

The maps $\Psi$ and $\Phi$ exhibit the following structural invariances, which reflect the geometry of the spaces on which they act:

\begin{enumerate}
\item \textbf{Row-wise translation invariance.}  
Adding a constant scalar to every entry of any row of $P$ leaves $\Psi(P)$ unchanged, because the centering operator $H$ removes all row means.

\item \textbf{Row-wise scaling invariance.}  
Multiplying any row of $P$ by a positive scalar does not affect $\Psi(P)$: centering eliminates the additive part, and subsequent normalization removes the multiplicative factor.

\item \textbf{Row sign symmetry.}  
Flipping the sign of any row of $\Psi(P)$ induces the corresponding sign flip in both the row and column of $\Phi(\Psi(P))$. Thus $\Phi$ is equivariant under diagonal sign matrices.

\item \textbf{Right-orthogonal invariance.}  
If $Z \in \Mn$ and $Q$ is any orthogonal matrix, then $\Phi(ZQ) = ZZ^\top$. In
particular, if $Q$ acts as an orthogonal change of basis on $\ones^\perp$
(i.e., $Q^\top \ones = \ones$), then $\Phi$ depends only on the row
configuration in $\ones^\perp$, not on the chosen orthonormal basis of that
hyperplane.
\end{enumerate}

These invariances imply that the dynamics $P_{k+1} = f(P_k)$ act naturally on the structured sets $\Mn$ and $\Sym$. Starting from a nondegenerate initialization $P_0$ drawn from \eqref{eq:init}, the iteration may be viewed as an alternating projection-like procedure in a space endowed with centering, normalization, and correlation symmetries.

\paragraph{Why these geometric preliminaries are needed.}
Although the present work is empirical, a minimal amount of geometric
structure is essential for interpreting the iteration and for ensuring that
all diagnostics are well defined.  
The map \(f=\Phi\circ\Psi\) requires that every centered row be nonzero and
that each iterate remain a correlation matrix with unit diagonal.  
Forward invariance of the elliptope $\mathcal{E}_n$, together with the standing 
nondegeneracy assumption \eqref{eq:nondeg-4}, guarantees precisely this: 
the iteration never leaves the region on which $\Psi$ and $\Phi$ are well defined, 
so quantities such as the Frobenius step $\Delta_k$, the entrywise step $E_k$, 
and the contraction ratio $\rho_k$ are comparable across iterations and dimensions.
 
Moreover, the local geometric behavior of \(f\)—strong shrinkage in directions
that alter row means or row norms, and near-isometry in directions that keep
row means equal to zero and row norms equal to one—provides the conceptual
foundation for the empirical laws presented in Section~\ref{sec:empirical}.  
These preliminaries therefore supply the structural context within which the
numerical results can be interpreted.

\section{Experimental Setup and Quantitative Diagnostics}
\label{sec:metrics-experiments}

This section specifies the quantitative framework used to formulate the
empirical regularities and laws observed in \cref{sec:empirical}. The setup enables comparability across dimensions, and supports statistical validation of
the observed laws. Historical precedents and qualitative expectations can be
traced to \cite{mcquitty1968multiple,breiger1975clustering,chen2002gap},
and the norms and diagnostics follow standard practice in fixed-point
matrix-iteration theory \cite{bauschke2017convex,sinkhorn1964relationship}.

\subsection{Initialization and degeneracy checks}\label{subsec:init}

For each dimension \(n\in\{3,\dots,2000\}\) and trial \(t=1,\dots,N\) with
\(N\in[10^2,10^3]\), an i.i.d.\ matrix
\(P_0^{(t)}\in[-1,1]^{n\times n}\) is drawn and the first step
\(P_1^{(t)}=f(P_0^{(t)})\) is computed. The row-sphericalization \(\Psi\)
requires nondegeneracy of centered rows:
\[
\qquad \|(P_0^{(t)}H)_{i,:}\|_2>0\quad\text{for all }i.
\]
Trials that fail this test have probability zero in the continuous model and
are numerically rare due to finite precision. Such trials are discarded and
resampled. This guarantees \(P_1^{(t)}\in\Sym\) and hence
\(P_k^{(t)}\in\Sym\) for all \(k\ge 1\), as described in subsection~\ref{subsec:forward}.

\subsection{Metrics and norm equivalence}\label{subsec:metrics}

At each iterate the Frobenius step size, the maximal entrywise step, and
the one-step contraction ratio are recorded:
\begin{equation}\label{eq:delta-E-rho}
\Delta_k:=\norm{P_{k+1}-P_k}_F,\qquad
E_k:=\max_{i,j}\lvert P_{k+1}(i,j)-P_k(i,j)\rvert,\qquad
\rho_k:=\frac{\Delta_{k+1}}{\Delta_k}\ \ (\Delta_k>0).
\end{equation}

By the norm-equivalence relation noted in Subsection~\ref{subsec:iteration-operator}, the sequences \((\Delta_k)\) and \((E_k)\) capture the same qualitative decay behavior up to the fixed factor $n$.

\subsection{Stopping rules and convergence counters}\label{subsec:stopping}

For a tolerance \(\varepsilon>0\), convergence is declared once
\[
\Delta_k<\varepsilon.
\]
The convergence time is defined as
\begin{equation}\label{eq:T-def}
T(n,\varepsilon):=\min\bigl\{k:\Delta_j<\varepsilon\ \text{for all } j\ge k\bigr\}.
\end{equation}
Results are reported for $\varepsilon \in \{10^{-6},10^{-12}\}$.
The tolerance $10^{-6}$ reflects standard numerical accuracy, whereas
$10^{-12}$ is chosen as a near–practical resolution limit for double‐precision
iterative computations: although machine epsilon is approximately
$2.22\times 10^{-16}$, accumulated rounding and cancellation typically
make $10^{-12}$ the deepest reliable threshold.  
This allows the iteration to be tracked down to the point where rounding
effects begin to dominate.

\subsection{Implementation and computational cost}
\label{subsec:implementation}

The empirical iteration applies the correlation map $P_{k+1}=f(P_k)$ by
explicitly evaluating the centered covariance between all pairs of rows.
Given $P_k\in\mathbb{R}^{n\times n}$, the code computes the mean and centered
variance of each row, forms the centered covariance between rows $i$ and $j$,
normalizes it, and assigns the value symmetrically to
$(P_{k+1})_{ij}=(P_{k+1})_{ji}$. The diagonal of $P_{k+1}$ is reset to one
after each update so that all iterates remain correlation matrices.

The computational cost of a single update scales cubically with $n$.
Computing row means and variances requires $O(n)$ work per row, and each of
the $\tfrac{n(n+1)}{2}$ centered covariances costs a further $O(n)$
operations, yielding a total cost of $O(n^{3})$ floating-point operations per iteration.  Storing two successive iterates together with a small set of working quantities requires $O(n^{2})$ memory.

The iteration terminates once
\[
\max_{i,j}\bigl|(P_{k+1})_{ij}-P_{k}(i,j)\bigr|\;\le\;\varepsilon,
\qquad \varepsilon = 10^{-12},
\]
or after a safeguard limit of $1000$ updates.

To obtain statistics across initial conditions, the iteration is repeated for
$N$ independently sampled matrices $P_0$, each using an independent random
seed.  For every run, the full diagnostic trace $(\Delta_k,E_k,\rho_k)$ and the
corresponding convergence time $T(n,\varepsilon)$ are recorded.

\paragraph{Breadth of initialization.}
The empirical behavior observed below persists across an extremely large
space of starting matrices. A symmetric $n\times n$ correlation matrix has
$\tfrac{n(n-1)}{2}$ free entries, so even a coarse grid with resolution
$10^{-3}$ already yields $2001^{\,n(n-1)/2}$ admissible initial states
(e.g.\ $\sim 10^{9}$ for $n=3$ and $\sim 10^{50}$ for $n=6$). The i.i.d.\
floating-point initialization used here explores an even larger region, yet
all trajectories exhibit the same qualitative structure across dimensions and
trials. This indicates that the empirical laws reported later are not tied to
a small or atypical initialization set.

\subsection{Reproducibility checklist}\label{subsec:repro}

\begin{itemize}
  \item \textbf{Input.}
        Dimensions $n\in\{3,\dots,2000\}$; trials per dimension $N\in[10^2,10^3]$; initial matrices $P_0$ sampled independently from \eqref{eq:init}.

  \item \textbf{Iteration.}
        Update rule $P_{k+1}=f(P_k)$ as in \eqref{eq:Pk-plus-one}.  
        Any $P_0$ whose centered rows contain a zero vector (a probability–zero event under \eqref{eq:init}) is discarded and resampled.

  \item \textbf{Diagnostics.}
        Per iteration, record $\Delta_k$, $E_k$, and $\rho_k$.  
        Convergence is declared when $\Delta_k<\varepsilon$ for $\varepsilon\in\{10^{-6},10^{-12}\}$, and $T(n,\varepsilon)$ is computed using \eqref{eq:T-def}.

  \item \textbf{Outputs.}
        For each $(n,t)$ pair, store the iterates $\{P_k\}$, the diagnostic traces $(\Delta_k,E_k,\rho_k)$, and the derived statistics (trajectory plots, pooled distributions, and summaries of $T(n,\varepsilon)$).  
        All figures, tables, and raw datasets are produced automatically by the accompanying scripts.
\end{itemize}

The full experimental framework—including implementation, execution scripts, and figure-generation utilities—is permanently archived on Zenodo \cite{alhajjhassan2025framework}, enabling independent reproduction of all reported results.


\section{Empirical Laws and Observed Regularities}\label{sec:empirical}

Using the metrics and diagnostics from \cref{sec:metrics-experiments}, four
robust empirical regularities of the discrete flow \(P_{k+1}=f(P_k)\) are
now described. These laws formalize and extend classical qualitative
observations in psychometrics, network block-modeling, and association plots
\cite{mcquitty1968multiple,breiger1975clustering,chen2002gap}, and clarify
the differences from contraction-based scalings such as Sinkhorn's
algorithm \cite{sinkhorn1964relationship}. For historical local results near
pattern matrices, see Kruskal’s notes~\cite{kruskalCONCOR,chen2002gap};
for analogies with iterative normalization in machine learning, see
\cite{huang2019iterative}.
\subsection{Law I: Universal first-step contraction}\label{subsec:law1}

\begin{empiricallaw}[Dimension-independent initial drop]\label{law:law1}
Let \(P_0 \in [-1,1]^{n\times n}\) have i.i.d.\ entries and no zero-centered
row, and let \(\Delta_k := \|P_{k+1}-P_k\|_F\). Then with high probability
\begin{equation}\label{eq:rho0-def}
\rho_0 \;:=\; \frac{\Delta_1}{\Delta_0} \;\ll\; 1 ,
\end{equation}
and the empirical distribution of \(\rho_0\) is essentially independent of
the dimension \(n\). Equivalently, for all tested sizes the map
\(P_0 \mapsto P_1\) exhibits a sharp first-step contraction in Frobenius
distance, with no loss of contraction as \(n\) increases.
\end{empiricallaw}

\begin{figure}[ht]
    \centering
    \includegraphics[width=\linewidth]{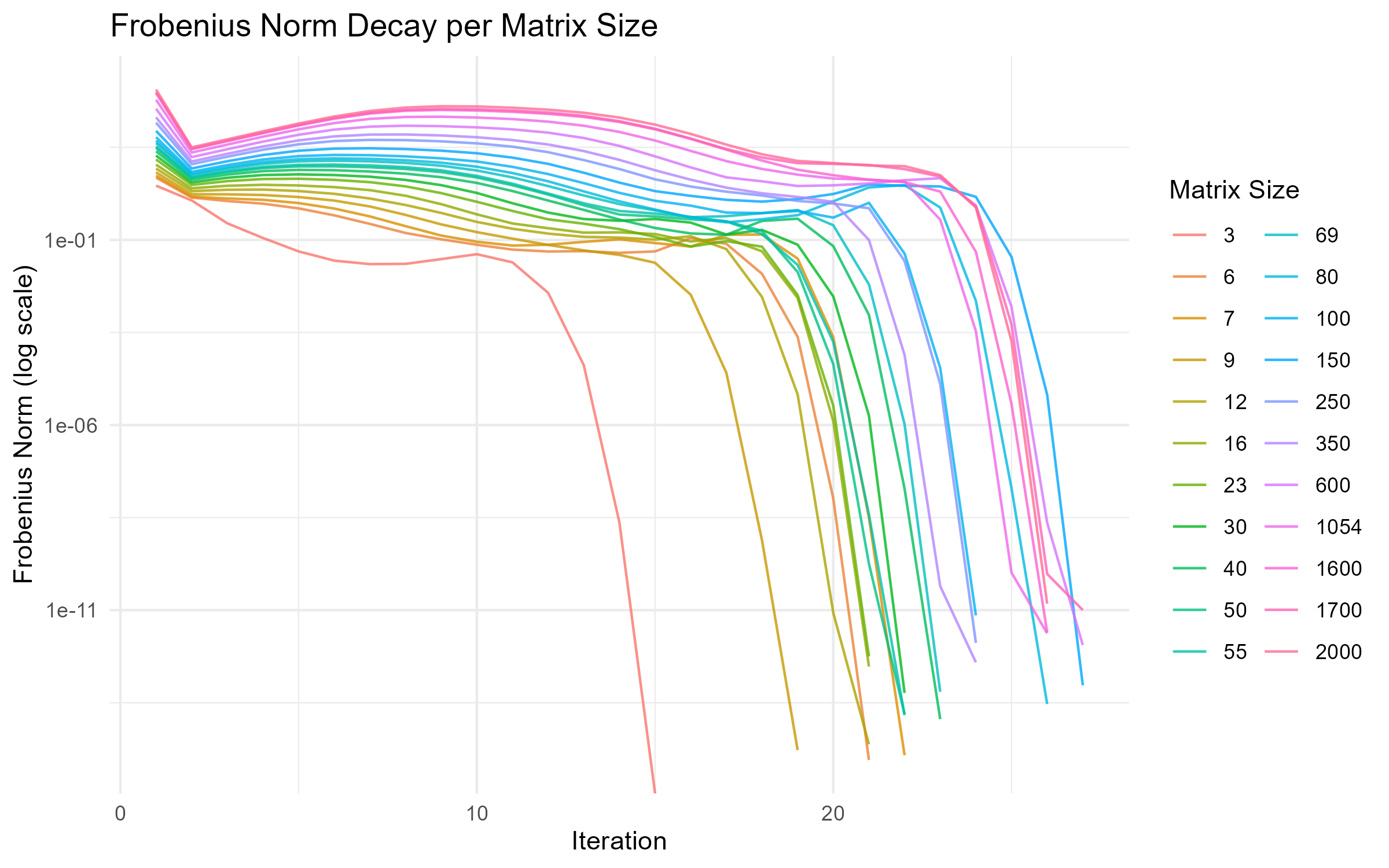}
    \caption{\textbf{Frobenius step-size decay across iterations.}
    Each curve corresponds to a fixed matrix dimension \(n\).
    All dimensions exhibit a pronounced contraction between
    \(k=0\) and \(k=1\), followed by a near-isometric regime.
    This universal pattern is the content of Law~I.}
    \label{fig:frob-decay-all}
\end{figure}

\begin{figure}[ht]
    \centering
    \includegraphics[width=\linewidth]{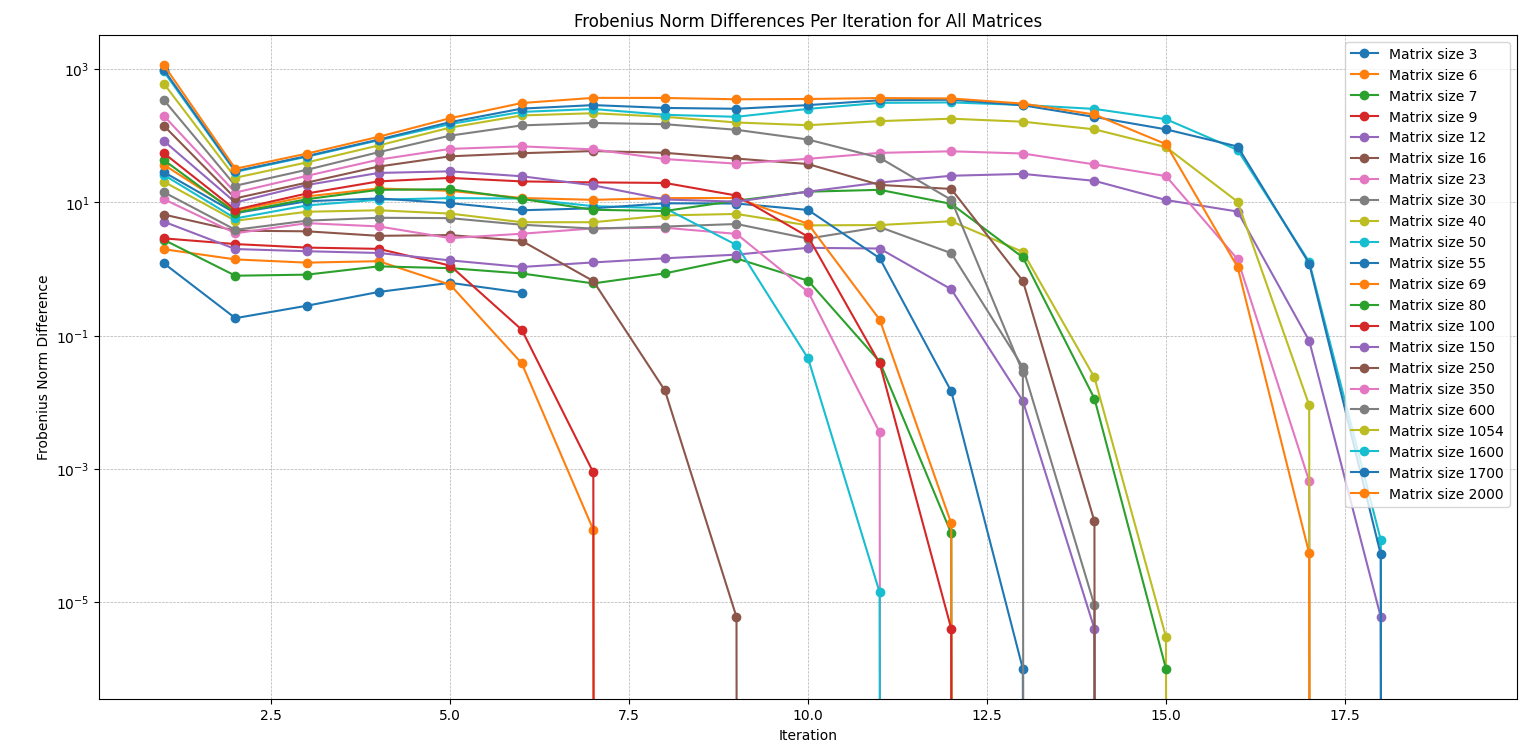}
    \caption{\textbf{Per iteration view of trajectories per matrix size.}
    A detailed view of iterations for multiple sizes \(n\).
    The Frobenius step size undergoes a large reduction at \(k=1\) for
    every dimension, illustrating that the sharp first-step contraction is
    uniform in \(n\).}
    \label{fig:frob-decay-highlight}
\end{figure}

\begin{figure}[ht]
  \centering
  \includegraphics[width=\linewidth]{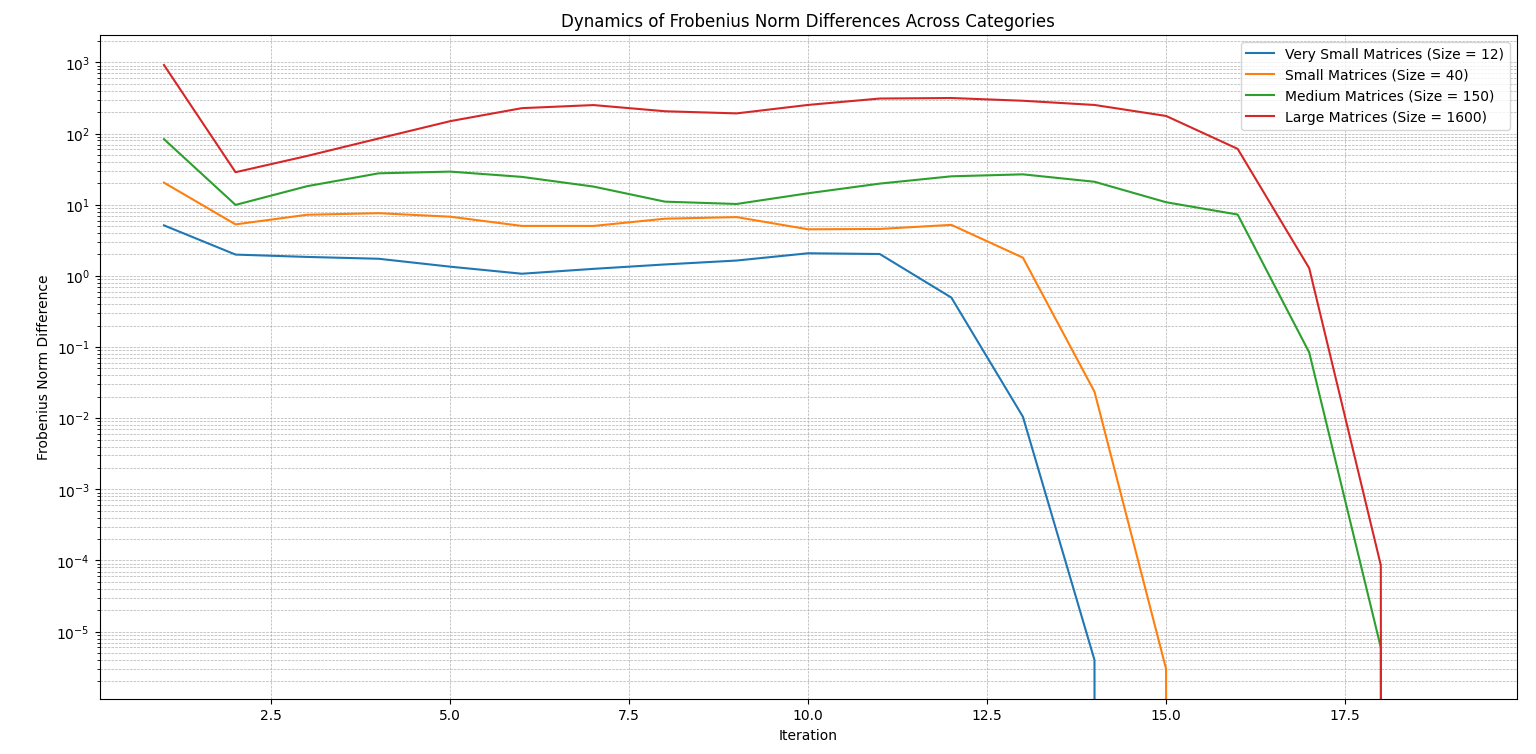}
  \caption{\textbf{Representative trajectories by matrix size.}
  Each trajectory shows the Frobenius step-size decay for a representative size from each dimension group. All four groups exhibit an immediate and sharp first-step
  contraction, confirming the dimension-independent nature of the
  phenomenon. Larger dimensions display a sharper drop and a quantitatively smaller ratio
  \(\rho_0\), but the qualitative behavior is identical across sizes.}
  \label{fig:frob-decay-by-size}
\end{figure}

\paragraph{Empirical evidence.}
Figure~\ref{fig:frob-decay-all} displays representative trajectories of the
Frobenius step size
\(\Delta_k = \|P_{k+1}-P_k\|_F\)\!,
for several dimensions \(n\). All trajectories share the same qualitative
behavior: a pronounced contraction at the first iteration
(\(k=0\to1\)), followed by a markedly slower decay over subsequent
iterations. Figure~\ref{fig:frob-decay-highlight} provides a magnified view of the initial region, highlighting the uniformity of the drop across sizes: the first step is sharply contractive for every tested \(n\in[3,2000]\). Figure~\ref{fig:frob-decay-by-size} groups trajectories by dimension category. The qualitative shape of the initial drop is stable across all sizes, whereas the precise magnitude of the ratio \(\rho_0\) decreases gradually with \(n\). The contraction phenomenon itself, however, occurs for all sizes without exception.

\paragraph{Numerical summary.}

\begin{figure}[ht]
  \centering
  \includegraphics[width=\textwidth]{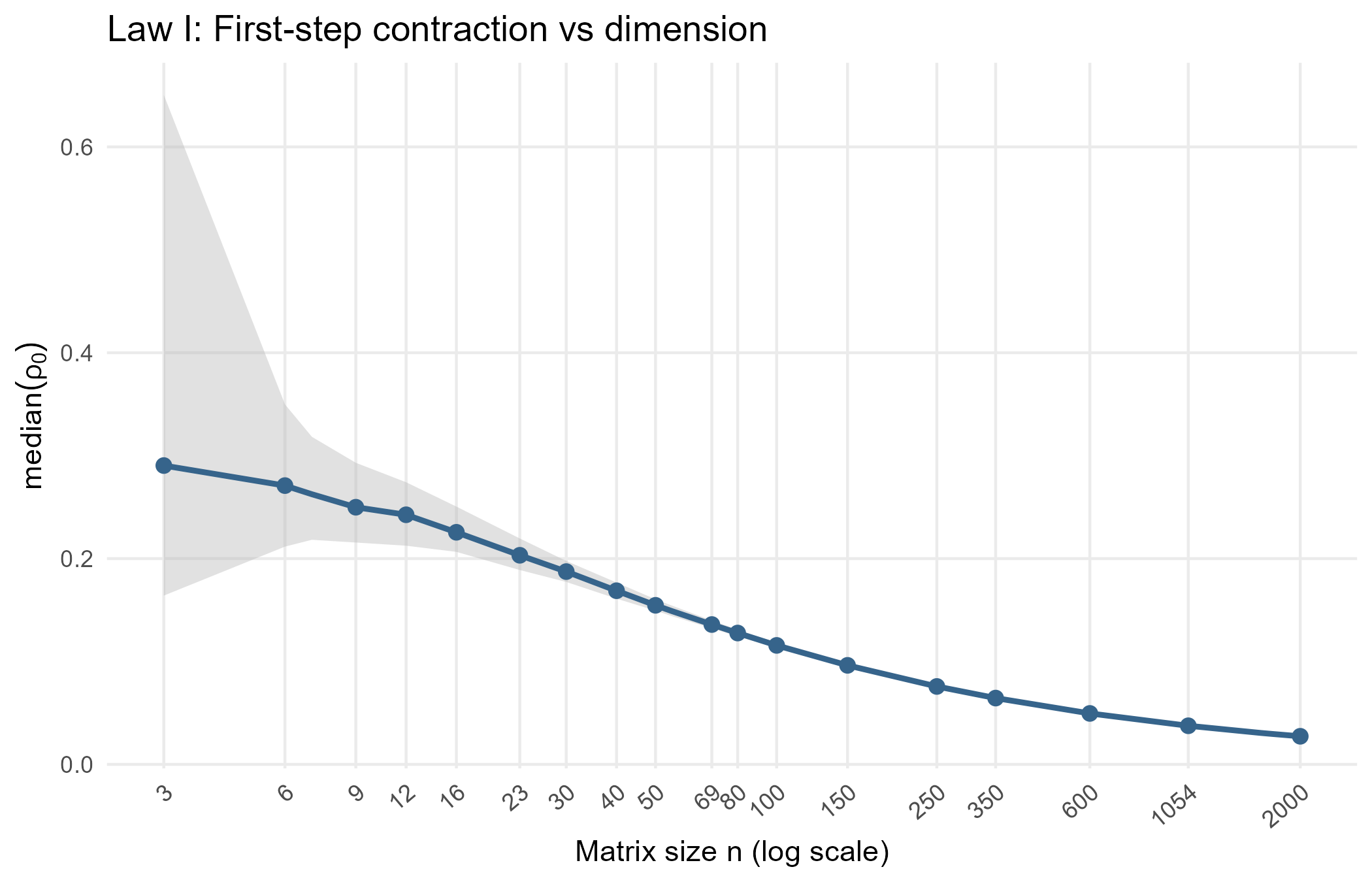}
  \caption{Median first-step contraction $\rho_0$ versus dimension $n$ (logarithmic horizontal axis).  
  Contraction occurs for every dimension and every trial, while its magnitude
  strengthens monotonically with $n$. The shaded region denotes the
  interquartile range, which decreases sharply for large $n$, indicating
  increasing stability of contraction strength.}
  \label{fig:law1-median-curve}
\end{figure}

\begin{table}[ht]
\centering
\small
\caption{Representative median first-step contraction ratios $\rho_0$ showing strengthening with dimension.}
\label{tab:law1-magnitude}
\begin{tabular}{r r r}
\toprule
$n$ & Median $\rho_0$ & IQR \\
\midrule
3   & 0.29 & [0.16,\,0.65] \\
12  & 0.24 & [0.21,\,0.27] \\
50  & 0.15 & [0.15,\,0.16] \\
150 & 0.10 & [0.09,\,0.10] \\
2000 & 0.03 & [0.03,\,0.03] \\
\bottomrule
\end{tabular}
\end{table}
Table~\ref{tab:law1-magnitude} reports representative values of the first-step
contraction ratio $\rho_{0}$, computed over $1000$ trials for each dimension.
In all $22\times 1000$ experiments we observed $\rho_{0}<1$, confirming that the
first iteration is universally contractive, independently of $n$.  At the same
time, the {\em strength} of this contraction increases systematically with
dimension: the median value of $\rho_{0}$ decreases by almost an order of
magnitude between $n=3$ and $n=2000$, and its interquartile range collapses to
nearly zero for large $n$.  The scaling behavior is clearly visible in
Figure~\ref{fig:law1-median-curve}, which depicts the median $\rho_{0}$ as a
function of $n$ on a logarithmic axis.  The data therefore show that initial
contraction is both dimension--independent in occurrence and increasingly
strong at higher dimensions.

\paragraph{Geometric explanation.}
The sharp and dimension--independent contraction observed at the first
iteration admits a transparent geometric explanation.  
Let \(P_0 \in [-1,1]^{n\times n}\) be the initial random matrix with
i.i.d.\ entries drawn from the uniform distribution on \([-1,1]\), and let
\(p_i \in \mathbb{R}^n\) denote its \(i\)-th row.
Since \(\mathbb{E}\,p_{ij}^2 = 1/3\), one has
\[
\mathbb{E}\,\|p_i\|_2^2 = \frac{n}{3},
\qquad\text{so}\qquad
\|p_i\|_2 = O(\sqrt{n})
\]
with high probability.  
Thus each row contains a large component in the direction of the all-ones
vector~\(\ones\).  The centering operator
\(H = I - \frac{1}{n}\ones\ones^{\mathsf T}\) removes this common component,
producing
\[
q_i := p_i - \bar p_i\,\ones ,
\qquad
\bar p_i := \frac{1}{n}\sum_{j=1}^n p_{ij}.
\]
A direct computation shows that
\[
\mathbb{E}\,\|q_i\|_2^2 = \frac{n-1}{3},
\]
so each residual vector \(q_i\) lies in the \((n-1)\)-dimensional hyperplane
\(\ones^\perp\) and still has norm \(O(\sqrt{n})\), but with its dominant
variance contribution along \(\ones\) removed.

A key point is that, by the law of large numbers, the row means \(\bar p_i\)
concentrate tightly around zero.  Consequently, the subtraction
\(p_i \mapsto q_i\) is highly coherent across~\(i\): each row loses almost the
same vector in the direction of \(\ones\), forcing the residuals to cluster
in \(\ones^\perp\) with near-isotropic orientations for large \(n\).

Normalizing these centered rows yields \(\Psi(P_0)\), whose rows are \(n\)
almost-independent unit vectors in an \((n-1)\)-dimensional space.  
In such a setting---placing \(n\) unit vectors into a space of dimension
\(n-1\)---generic pairs have inner products of order \(O(1/\sqrt{n})\).  
Hence, the correlation matrix
\[
P_1 = \Psi(P_0)\,\Psi(P_0)^{\mathsf T}
\]
satisfies, for \(i \neq j\),
\[
(P_1)_{ij} = O\!\left(\frac{1}{\sqrt{n}}\right).
\]

This dichotomy of scales explains the first-step contraction.
The transition \(P_0 \mapsto P_1\) removes a coherent component of magnitude
\(O(\sqrt{n})\) from each of the \(n\) rows, giving
\[
\Delta_0 = \|P_1 - P_0\|_F = O(n).
\]
By contrast, the next update compares two correlation matrices whose
off--diagonal entries already have magnitude \(O(1/\sqrt{n})\), and therefore
involves adjustments only of size \(O(\sqrt{n})\):
\[
\Delta_1 = O(\sqrt{n}).
\]
Thus the first-step contraction ratio satisfies
\begin{equation}\label{eq:rho0-scaling}
\rho_0
   \;=\;
   \frac{\Delta_1}{\Delta_0}
   \;=\;
   O\!\left(\frac{1}{\sqrt{n}}\right),
\end{equation}
in agreement with the empirical decay of the medians reported in
Table~\ref{tab:law1-magnitude}.  
This high--dimensional geometric structure explains why the \emph{shape} of
the first-step contraction is independent of dimension, while its
\emph{magnitude} becomes more pronounced as \(n\) increases.

\paragraph{Conclusion.}
Across all \(22\,000\!+\) trials and all tested dimensions, the first-step
ratio \(\rho_0\) remained strictly below \(1\) and was typically far smaller.
The qualitative phenomenon of a sharp initial contraction occurs for every
dimension \(n\), while the quantitative magnitude of the contraction follows
the high-dimensional scaling estimate~\eqref{eq:rho0-scaling}. The first iteration
thus places the system very close to the regime in which the subsequent
near-isometric phase begins, and it does so uniformly for all matrix sizes.

\subsection{Law II: Nearly monotone decay, bounded overshoots, and controlled total variation}\label{subsec:law2}

\begin{empiricallaw}[Nearly monotone convergence with bounded overshoots]\label{law:law2}
For \(k\ge 1\), the ratios \(\rho_k\) satisfy \(\rho_k\approx 1\) with
bounded overshoots, and the overshoot bound vanishes as \(\Delta_k\to 0\).
More precisely, there exists a sequence \(\eta_k\downarrow 0\) such that
typically
\[
\Delta_{k+1}\le (1+\eta_k)\,\Delta_k,
\]
with occasional values \(\rho_k>1\), but with
\(\sup_k \rho_k<1+\varepsilon_\star(n)\), where \(\varepsilon_\star(n)\) is
small and does not grow with \(n\).
\end{empiricallaw}

For each trial \(t=1,\dots,N\) and iteration \(k\ge 0\), define the
Frobenius step
\begin{equation}\label{eq:frob-step}
\Delta_k^{(t)} \;=\; \bigl\| P_{k+1}^{(t)} - P_k^{(t)} \bigr\|_F,
\qquad
\bar{\Delta}_k \;=\; \frac{1}{N}\sum_{t=1}^N \Delta_k^{(t)} ,
\end{equation}
and the elementwise step
\begin{equation}\label{eq:elem-step}
E_k^{(t)} \;=\; \max_{i,j}\bigl\lvert P_{k+1}^{(t)}(i,j) - P_k^{(t)}(i,j) \bigr\rvert,
\qquad
\bar{E}_k \;=\; \frac{1}{N}\sum_{t=1}^N E_k^{(t)} .
\end{equation}
By \cref{eq:norm-equivalence}, the inequalities
\(E_k \le \Delta_k \le n\,E_k\) hold, so both metrics encode the same
qualitative decay up to constants.

\begin{figure}[htbp]
  \centering
  \includegraphics[width=\textwidth]{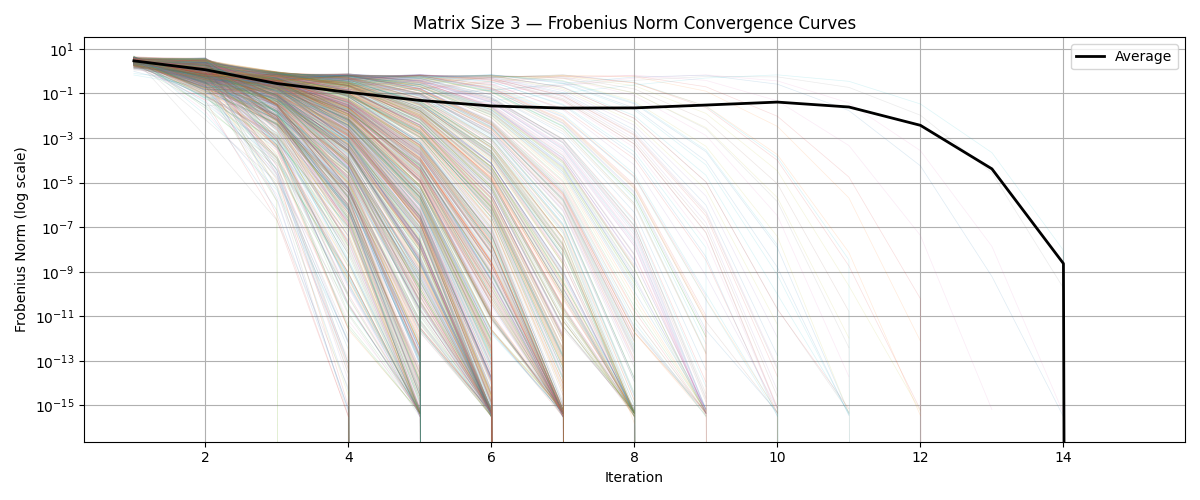}
\caption{\textbf{Frobenius step decay at $n=3$ over $N=1000$ trials.}
Colored curves show individual trajectories $\{\Delta_k^{(t)}\}$.  
The black curve is the mean $\bar\Delta_k$ over trials that are still 
active at iteration $k$.  
A universal sharp drop is visible at $k=1$, followed by near-monotone decay 
with shrinking oscillations.  
The mean curve ends at $k=14$, indicating that all trajectories have 
converged by this iteration.}
\label{fig:frob-n3-law2}
\end{figure}

\begin{figure}[htbp]
  \centering
  \includegraphics[width=\textwidth]{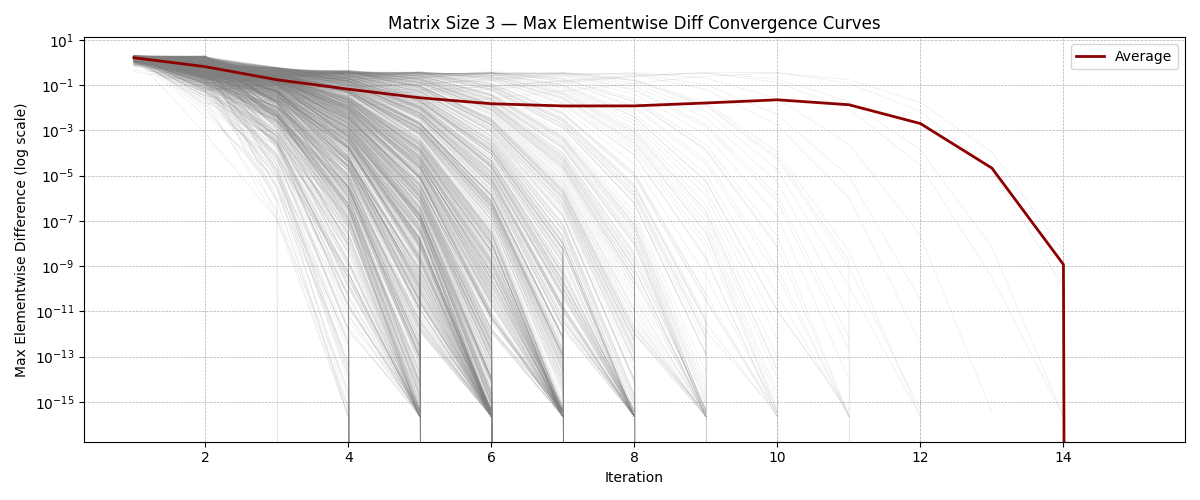}
 \caption{\textbf{Entrywise step decay at $n=3$ over $N=1000$ trials.}
Grey curves represent per-trial maximum entrywise difference $E_k^{(t)}$.  
The red curve is the evolving mean $\bar{E}_k$ over active trajectories.  
Its behavior mirrors that of $\Delta_k$: one large initial contraction, 
followed by smooth decay. The mean becomes invisible beyond $k=14$, indicating that all runs have stabilized by that point.}
  \label{fig:elem-n3-law2}
\end{figure}

\vspace{0.9em}
A representative example of this behavior is shown at $n=3$ in
Figures~\ref{fig:frob-n3-law2}–\ref{fig:elem-n3-law2}.  
Across $N=1000$ trials both diagnostics exhibit the same profile: a universal
first contraction, followed by smooth decay with small oscillations that
shrink as $k$ increases.  
The final disappearance of the mean curves at $k\approx 14$ merely reflects
that all trajectories have stabilized by this iteration; this observation
will be revisited later under Law~4, which concerns the boundedness of the
total convergence time.

\vspace{3em}

\vspace{6em}

\vspace{0.9em}

\begin{figure}[htbp]
  \centering
  \begin{subfigure}{0.32\linewidth}
    \centering
    \includegraphics[width=\linewidth]{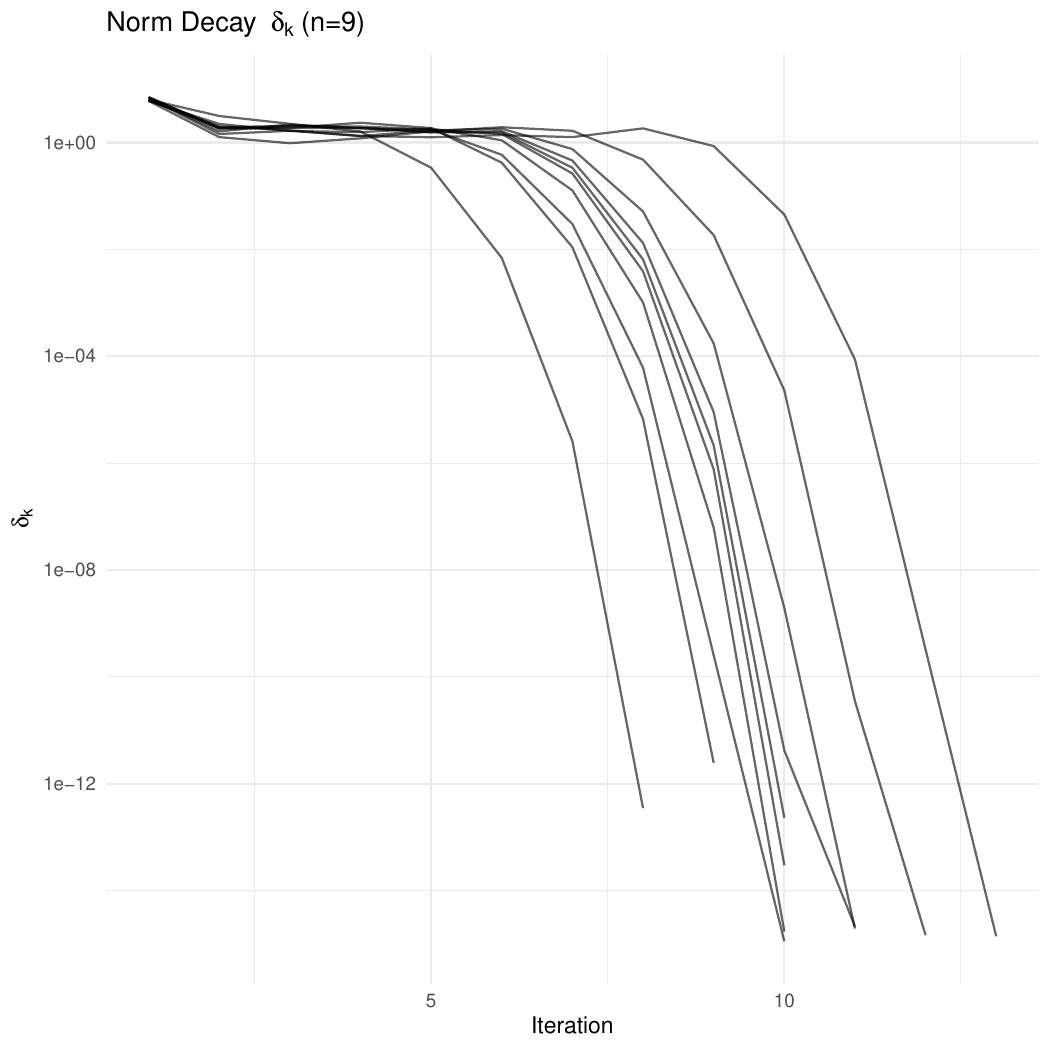}
    \caption{$n=9$}
  \end{subfigure}\hfill
  \begin{subfigure}{0.32\linewidth}
    \centering
    \includegraphics[width=\linewidth]{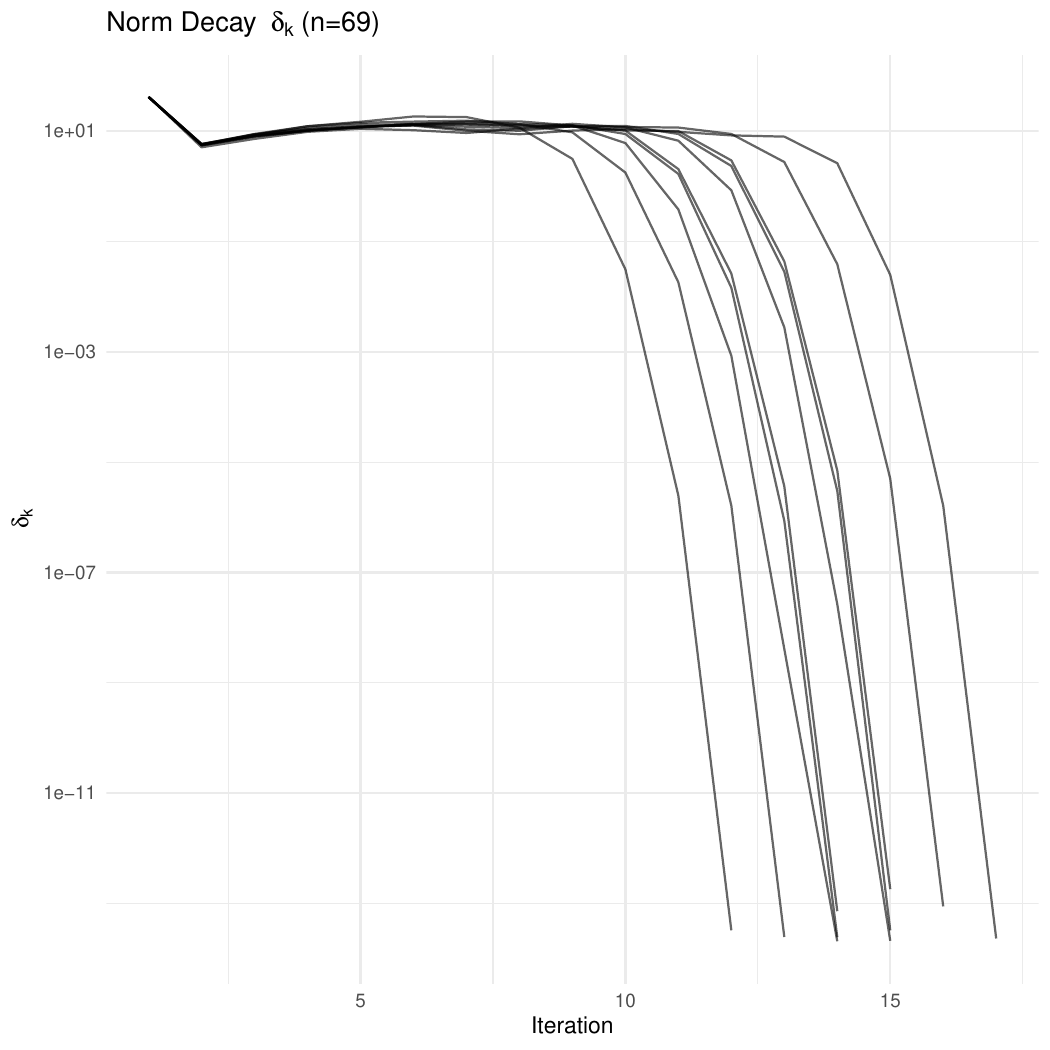}
    \caption{$n=69$}
  \end{subfigure}\hfill
  \begin{subfigure}{0.32\linewidth}
    \centering
    \includegraphics[width=\linewidth]{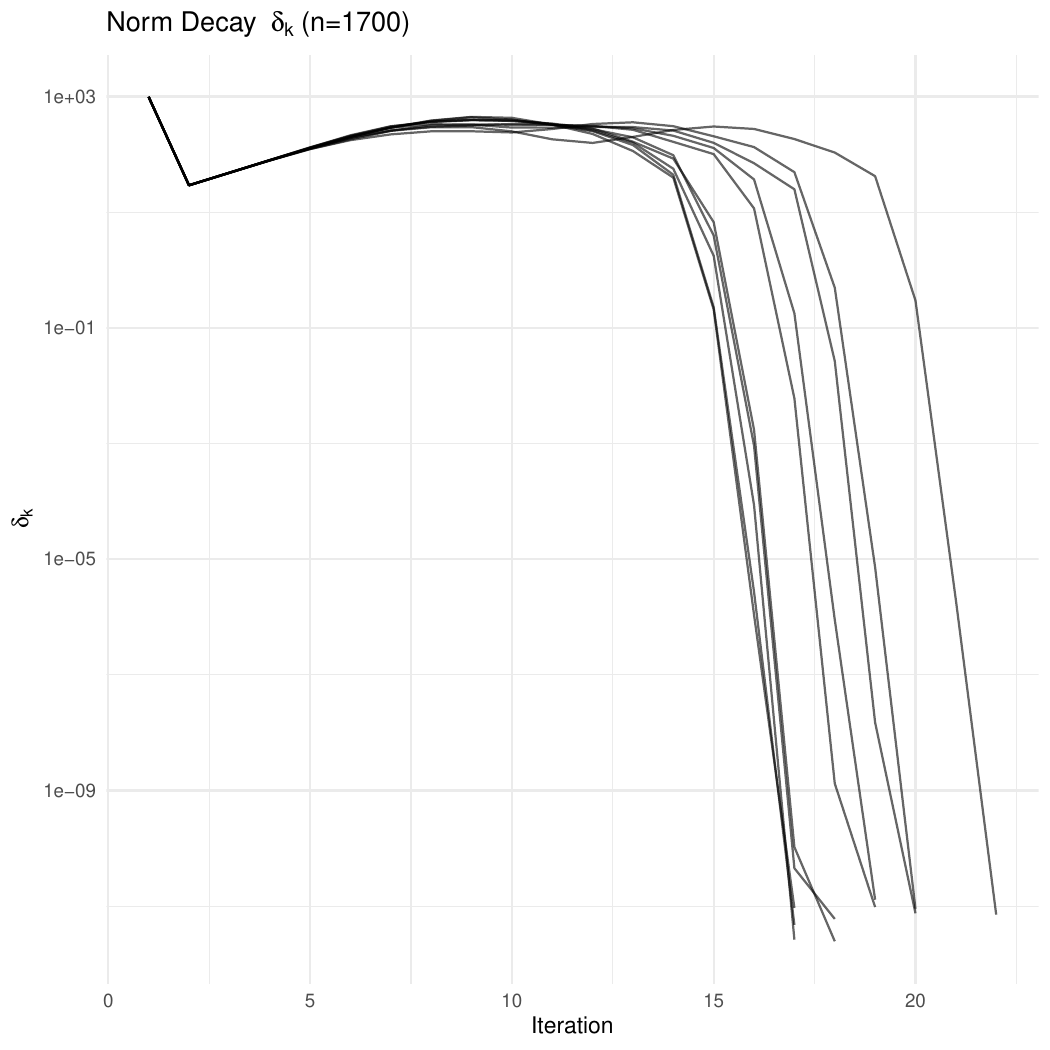}
    \caption{$n=1700$}
  \end{subfigure}
  \caption{\textbf{Frobenius step-size decay \(\Delta_k\) as a function of \(k\).}
  The trials across sizes display the same pattern: a universal first-step
  drop followed by nearly monotone decay with small bounded oscillations.
  The oscillation bound tightens with \(k\), consistent with
  \(\eta_k\downarrow 0\) in \cref{law:law2}.}
  \label{fig:norm-decay-examples-law2}
\end{figure}
\begin{figure}[htbp]
  \centering
  \includegraphics[width=0.90\linewidth]{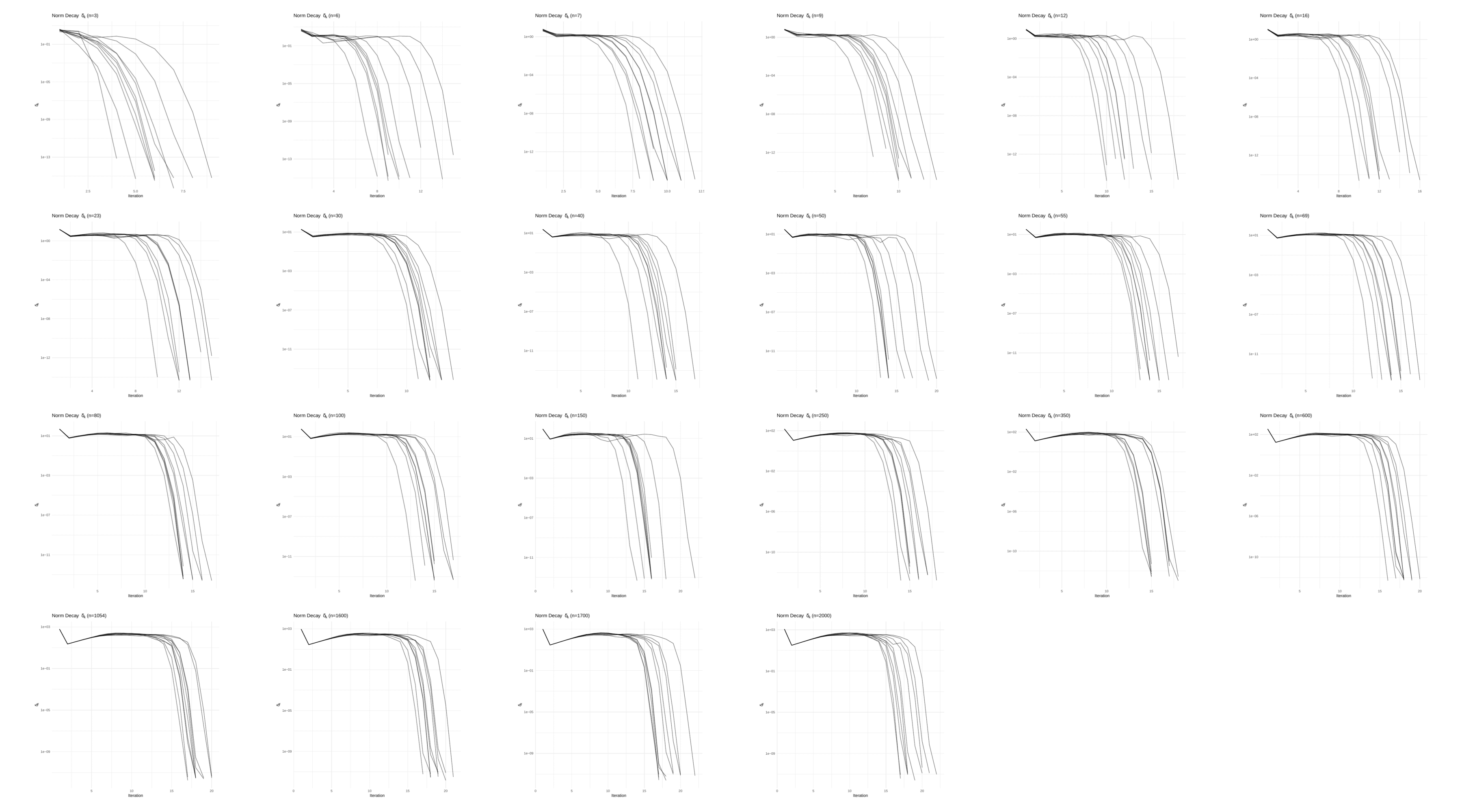}
  \caption{\textbf{Representative decay of \(\Delta_k\) across matrix sizes.}
For each dimension, a few representative trajectories are shown.
Across the full set of 1000 trials per dimension, the curves consistently
collapse to a common shape: a large initial descent followed by a
quasi-monotone phase with bounded overshoots. This behavior confirms
Law~II across the full grid of matrix sizes.}
\label{fig:norm-decay-grid-law2}
\end{figure}

\vspace{0.9em}
Representative Frobenius–norm decay trajectories for $n=3$, $n=69$, and $n=1700$, each sampled from $N=1000$ trials, are shown in \Cref{fig:norm-decay-examples-law2}.  
A full grid of trajectories across all tested dimensions is displayed in \Cref{fig:norm-decay-grid-law2}.  
In every case, the curves exhibit the same characteristic geometry: a strong universal first contraction followed by nearly monotone decay with small, uniformly bounded oscillations.
\vspace{0.9em}

\FloatBarrier

\begin{figure}[htbp]
  \centering
  \begin{subfigure}{0.48\linewidth}
    \centering
    \includegraphics[width=\linewidth]{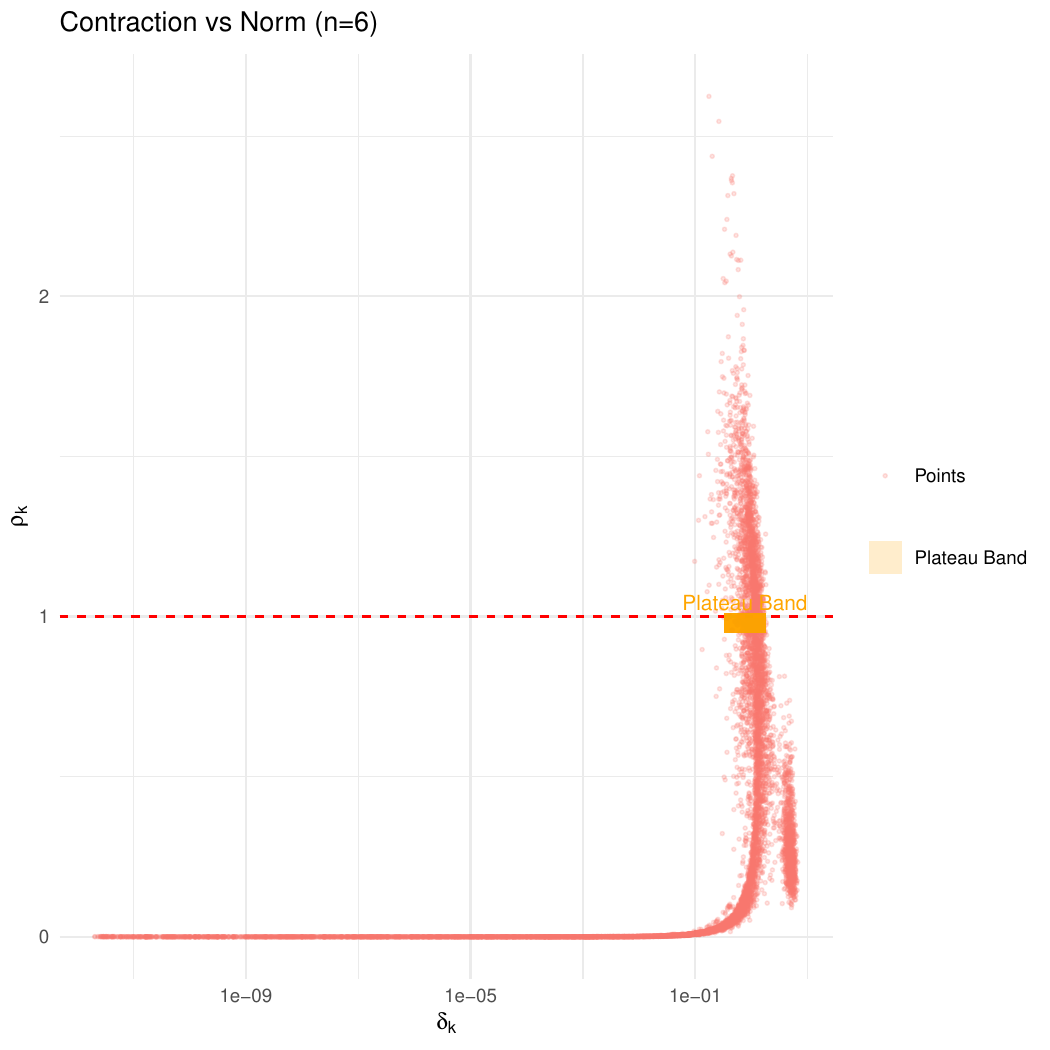}
    \caption{$n=6$}
  \end{subfigure}\hfill
  \begin{subfigure}{0.48\linewidth}
    \centering
    \includegraphics[width=\linewidth]{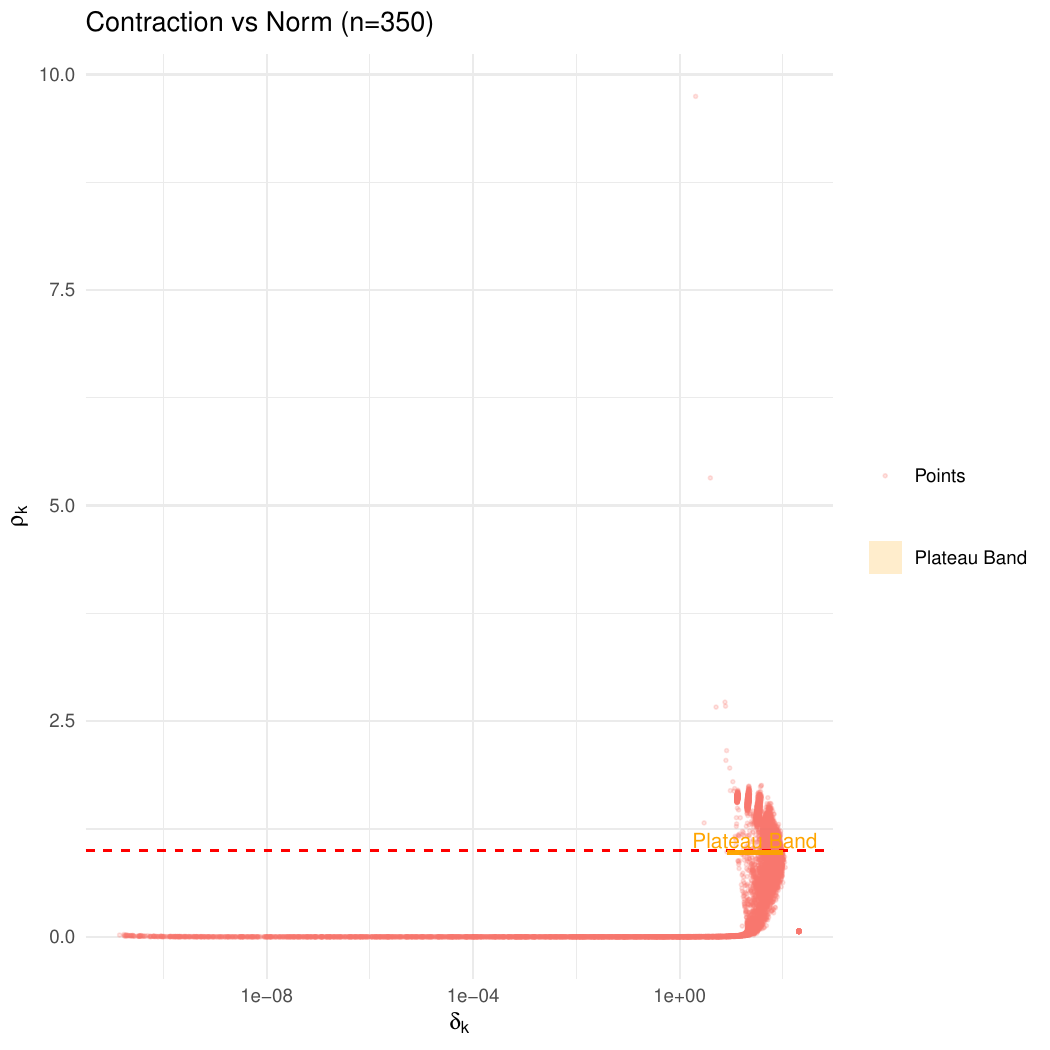}
    \caption{$n=350$}
  \end{subfigure}
    \caption{\textbf{Scatter of \(\rho_k\) versus \(\Delta_k\) at two sizes.}
  Each point corresponds to a single iterate \(k\ge 1\) from one trial. The
  horizontal axis is the step size \(\Delta_k=\norm{P_{k+1}-P_k}_F\); the
  vertical axis is the contraction ratio \(\rho_k=\Delta_{k+1}/\Delta_k\).
  For large \(\Delta_k\) in the early iterates, \(\rho_k\ll 1\), indicating
  strong contraction. As \(\Delta_k\downarrow 0\) in later iterates,
  \(\rho_k\to 1\) and the points remain confined in a narrow band above and
  below \(1\). This matches the quasi-monotone decay with bounded
  oscillations asserted in \cref{law:law2}.}
  \label{fig:rho-vs-delta-law2}
\end{figure}

\bigskip
\begin{figure}[htbp]
  \centering
\includegraphics[width=0.9\linewidth]{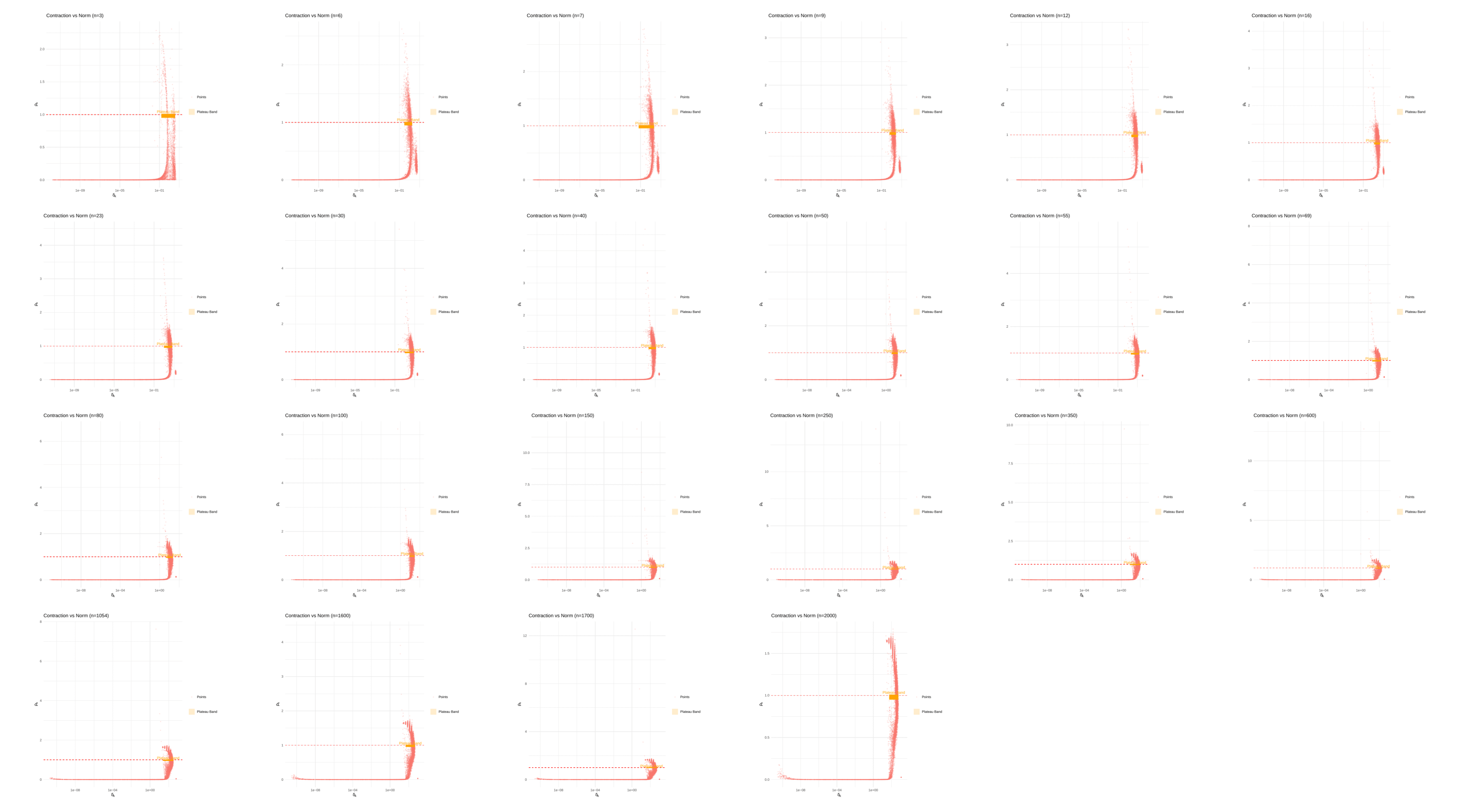}
\caption{\textbf{Representative scatter plots of \(\rho_k\) versus \(\Delta_k\) across matrix sizes.}
The grid displays the scatter plots of \((\Delta_k,\rho_k)\) for all tested
dimensions. For every size \(n\), the points exhibit the same distinctive
pattern: large values of \(\Delta_k\) correspond to strong contraction
(\(\rho_k \ll 1\)), while as \(\Delta_k \to 0\) the values of \(\rho_k\)
concentrate tightly near \(1\). The persistence of this structure across all
dimensions shows that the overshoot behavior is essentially independent of
matrix size.}
\label{fig:rho-vs-delta-grid-law2}
\end{figure}

\begin{figure}[htbp]
  \centering
  \includegraphics[width=0.9\linewidth]{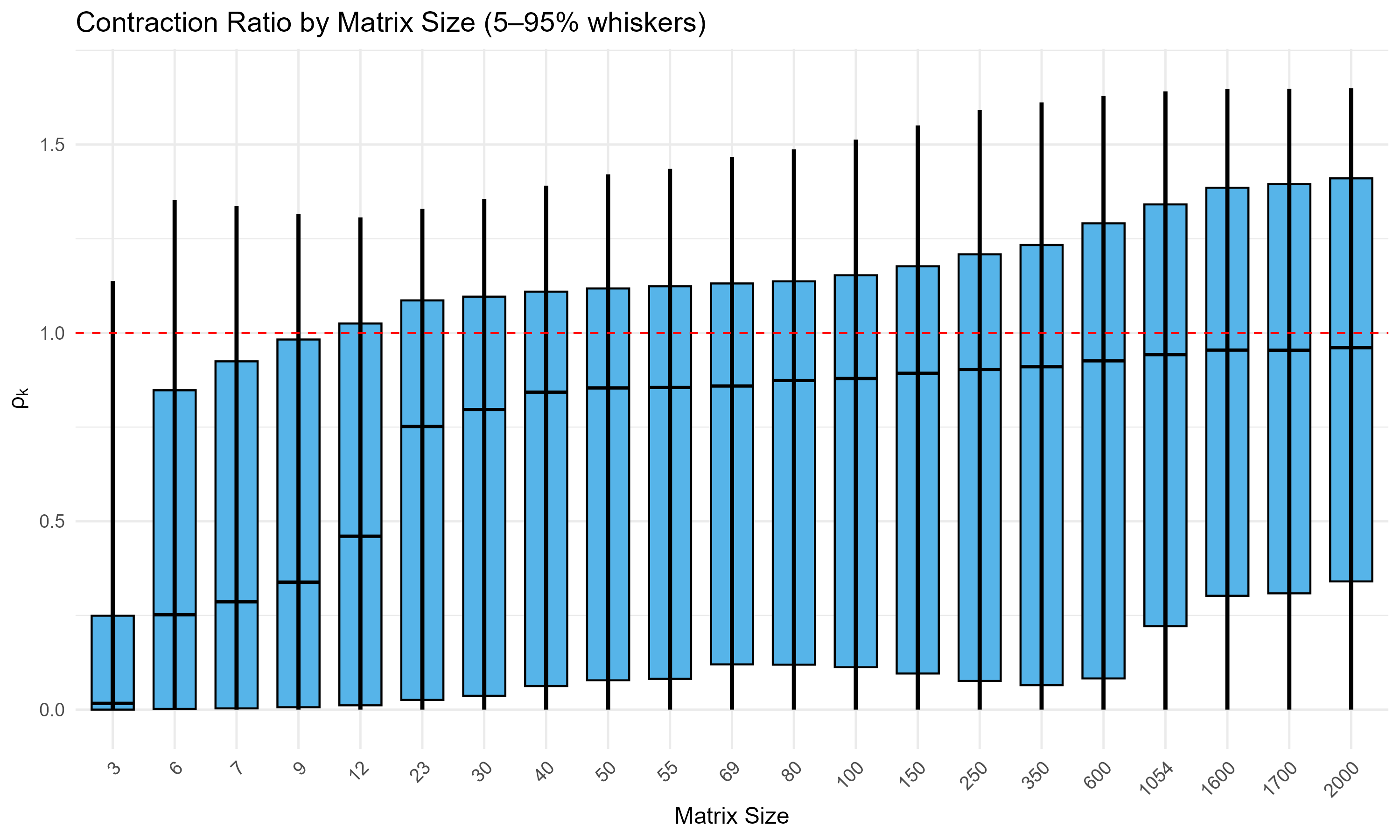}
\caption{\textbf{Contraction ratios \(\rho_k=\Delta_{k+1}/\Delta_k\) across matrix sizes (5--95\% whiskers).}
For each dimension \(n\), the box shows the interquartile range 
\([Q_{0.25},Q_{0.75}]\), with the thick bar indicating the median \(Q_{0.50}\); 
vertical whiskers correspond to the \(Q_{0.05}\) and \(Q_{0.95}\) quantiles. 
The dashed line \(y=1\) marks perfect isometry. Medians lie close to \(1\), 
indicating that, after the first iteration, most steps are non-expansive or 
nearly isometric. The interquartile ranges remain uniformly narrow across all 
sizes, and the 95\% whiskers stay close to \(1\), showing that the overshoot 
behavior is consistently small and largely independent of dimension.}

\label{fig:rho-boxplot-law2}
\end{figure}

Law~II is weaker than Fejér or quasi-Fejér monotonicity, because small
overshoots are observed, yet it remains compatible with the general
behavior of nonexpansive dynamics in Hilbert spaces~\cite{bauschke2017convex}.
Since the iteration is nonlinear and not uniformly contractive, classical
global convergence theorems for nonexpansive mappings
\cite{goebel1990topics,opial1967weak} do not apply directly.  The map exhibits
strong contraction in directions that alter row means or row norms, while
behaving nearly isometrically along directions that preserve them.  This
anisotropic behavior limits the applicability of standard fixed-point
tools.

\vspace{0.5em}
The bounded overshoots documented above, together with the clear tendency of $\rho_k$ to cluster tightly around $1$ as $\Delta_k \to 0$ (visible in Figures~\ref{fig:rho-vs-delta-grid-law2} and~\ref{fig:rho-boxplot-law2}), imply that the total variation along every trajectory
\begin{equation}
\label{eq:total-variation}
V_\infty \;=\; \sum_{k \ge 0} \Delta_k,
\end{equation}
is finite in all \(22{,}000+\) trials.

\vspace{0.5em}

Consistently with the heuristic scaling \(\Delta_0 = O(n)\) in Law~I and the uniformly bounded number of subsequent iterations in Law~IV, the total variation grows approximately linearly with the dimension: its median is about \(2.3\,n\) in Frobenius norm across the tested range (from roughly \(4.5\) for \(n=3\) to about \(4.6\times 10^3\) for \(n=2000\); see Table~\ref{tab:total-variation-summary}).

\vspace{0.5em}
\begin{table}[htbp]
\centering
\caption{Empirical distribution of the total variation \(V_\infty = \sum_{k\ge 0} \Delta_k\) (Frobenius norm) across all trials for selected dimensions.}
\label{tab:total-variation-summary}
\begin{tabular}{lrrrr}
\toprule
$n$     & Median & 95th Percentile & 99th Percentile & Max \\
\midrule
3       & 4.5    & 6.8        & 9.3        & 9.3 \\
100     & 217    & 256        & 286        & 312 \\
600     & 1360   & 1540       & 1660       & 1801 \\
2000    & 4580   & 5220       & 5580       & 6020 \\
\bottomrule
\end{tabular}
\end{table}
\vspace{0.5em}
\noindent This linear growth — dominated by the initial contraction — is fully consistent with the theoretical scaling from Law~I and the bounded iteration count that will be explored with Law~IV.

\begin{figure}[htbp]
  \centering
  \includegraphics[width=0.98\linewidth]{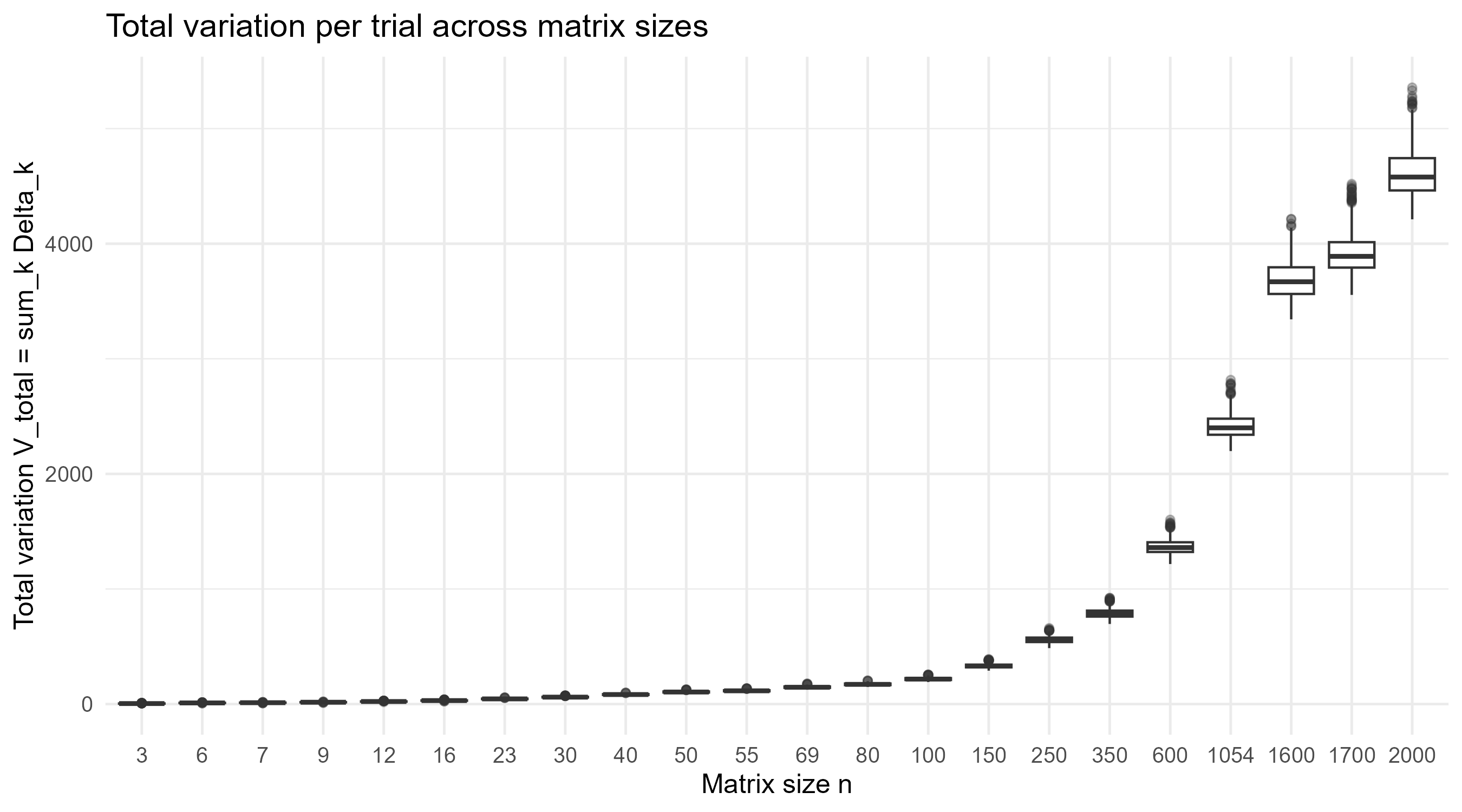}
  \caption{\textbf{Total variation \(V_\infty^{(t)}(n) = \sum_k \Delta_k^{(t)}\) by matrix size.}
  Boxplots over \(N=1000\) trials per size. The total variation is finite in
  every trajectory and grows roughly linearly with \(n\). Interquartile
  ranges remain narrow relative to the median at each dimension, indicating
  stable behavior across trials.}
  \label{fig:total-variation-boxplot}
\end{figure}

\paragraph{Total variation along trajectories.}
For each dimension \(n\) and trial \(t\), the cumulative change along the
trajectory is measured by the total variation
\[
V_\infty^{(t)}(n) \;=\; \sum_{k\ge 0} \Delta_k^{(t)}.
\]
Figure~\ref{fig:total-variation-boxplot} displays boxplots of
\(V_\infty^{(t)}(n)\) over \(N=1000\) trials for every tested size.
The total variation is finite in every run and grows approximately linearly
with \(n\): the median \(V_\infty^{(t)}(n)\) increases from about
\(4.5\) at \(n=3\) to about \(4.6\times10^3\) at \(n=2000\), with relatively
narrow interquartile bands around each median and maxima below \(5.4\times
10^3\) across all trials. Thus, although the cumulative motion is
dimension-dependent, it remains quantitatively controlled and exhibits no
sign of divergence.
\vspace{0.5em}

Taken together, these observations show that the iteration settles into a remarkably stable regime: oscillations remain uniformly bounded across all tested dimensions, and the values of $\rho_k$ concentrate in an increasingly tight band as $\Delta_k \to 0$. Decay progresses in a highly regulated fashion, with no runaway growth or irregular bursts—only a controlled descent that is effectively dimension-independent. This degree of regularity is striking for a nonlinear map with no \emph{a priori} monotonicity guarantees, and it reveals a structural coherence that any future analytical theory must account for.

\subsection{Law III: Dimension-independent contraction law}\label{subsec:law3}

The third empirical law describes how the instantaneous contraction ratio
\(\rho_k := \Delta_{k+1}/\Delta_k\) depends on the step magnitude
\(\Delta_k\). Across all tested matrix sizes, trials, and random
initializations, the cloud of points \((\Delta_k,\rho_k)\) exhibits the same
characteristic \emph{V-shaped} organization in the \((\Delta_k,\rho_k)\)-plane.
This structure, together with its overshoots, is an empirical finding of the
numerical study.

\paragraph{Binned median contraction curves.}
For each matrix size \(n\), all observed pairs \((\Delta_k,\rho_k)\) from all
trials are collected. The values of \(\Delta_k\) are partitioned into
logarithmic bins \(\{\mathcal B_b\}\). Within each bin we compute the median
step magnitude and the corresponding median contraction ratio
\begin{equation}\label{eq:rho-bar}
\delta_b := \operatorname{median}\{\Delta_k : \Delta_k \in \mathcal B_b\},
\qquad
\bar\rho_n(\delta_b)
   := \operatorname{median}\{\rho_k : \Delta_k \in \mathcal B_b\},
\end{equation}
together with the 10th and 90th percentiles of \(\rho_k\) in that bin.
The map \(\delta_b \mapsto \bar\rho_n(\delta_b)\) summarizes, for a fixed
dimension \(n\), how the contraction ratio varies as a function of the
current step size. The 10th–90th percentile band records the typical spread
of \(\rho_k\) around its median, including overshoots \(\rho_k>1\).

\begin{empiricallaw}[Universal dimension-independent contraction]\label{law:law3}
Across all dimensions \(n \in \{3,\dots,2000\}\), the curves
\(\delta_b \mapsto \bar\rho_n(\delta_b)\) and their 10th–90th percentile
bands exhibit the following empirical features:
\begin{enumerate}
  \item \textbf{Strong contraction for large steps.}
        For large \(\delta_b\) (corresponding to early iterations),
        the median contraction ratios satisfy \(\bar\rho_n(\delta_b)\ll 1\),
        and the 10th–90th percentile band lies well below~\(1\). Large steps
        are typically followed by much smaller steps.
  \item \textbf{Near-isometry with small overshoots for small steps.}
        As \(\delta_b \downarrow 0\), one observes
        \[
          \bar\rho_n(\delta_b) \longrightarrow 1,
        \]
        and the 10th–90th percentile band becomes a narrow strip around
        \(1\). Both \(\rho_k<1\) and \(\rho_k>1\) occur, but overshoots
        \(\rho_k>1\) remain small: their magnitude decreases with \(\delta_b\),
        and the envelope width shrinks as \(\delta_b \downarrow 0\).
  \item \textbf{Dimension-independence and V-shape.}
        When the curves \(\delta_b \mapsto \bar\rho_n(\delta_b)\) are plotted
        for all tested sizes, they lie almost exactly on top of one another.
        The combined graph forms a V-shaped pattern: a contractive branch
        with \(\bar\rho_n(\delta_b)\ll 1\) at large \(\delta_b\), and a
        near-horizontal, near-isometric branch with \(\bar\rho_n(\delta_b)
        \approx 1\) and small overshoots at small \(\delta_b\), with no
        visible dependence on \(n\).
\end{enumerate}
\end{empiricallaw}

\paragraph{Overshoots and spread.}
The 10th–90th percentile band around \(\bar\rho_n(\delta_b)\) provides a
quantitative description of overshoots. For intermediate and small step
sizes, the band typically straddles the line \(\rho=1\), indicating that
both slight contractions and slight expansions occur from one step to the
next. As \(\delta_b \downarrow 0\), this band contracts to a narrow interval
around \(1\), so overshoots \(\rho_k>1\) are present but uniformly small.
When all sizes are pooled before binning, the same pattern is observed:
the median curve and its envelope depend primarily on \(\Delta_k\), not on
the matrix dimension.

\begin{figure}[!htbp]
  \centering
  \includegraphics[width=0.92\linewidth]{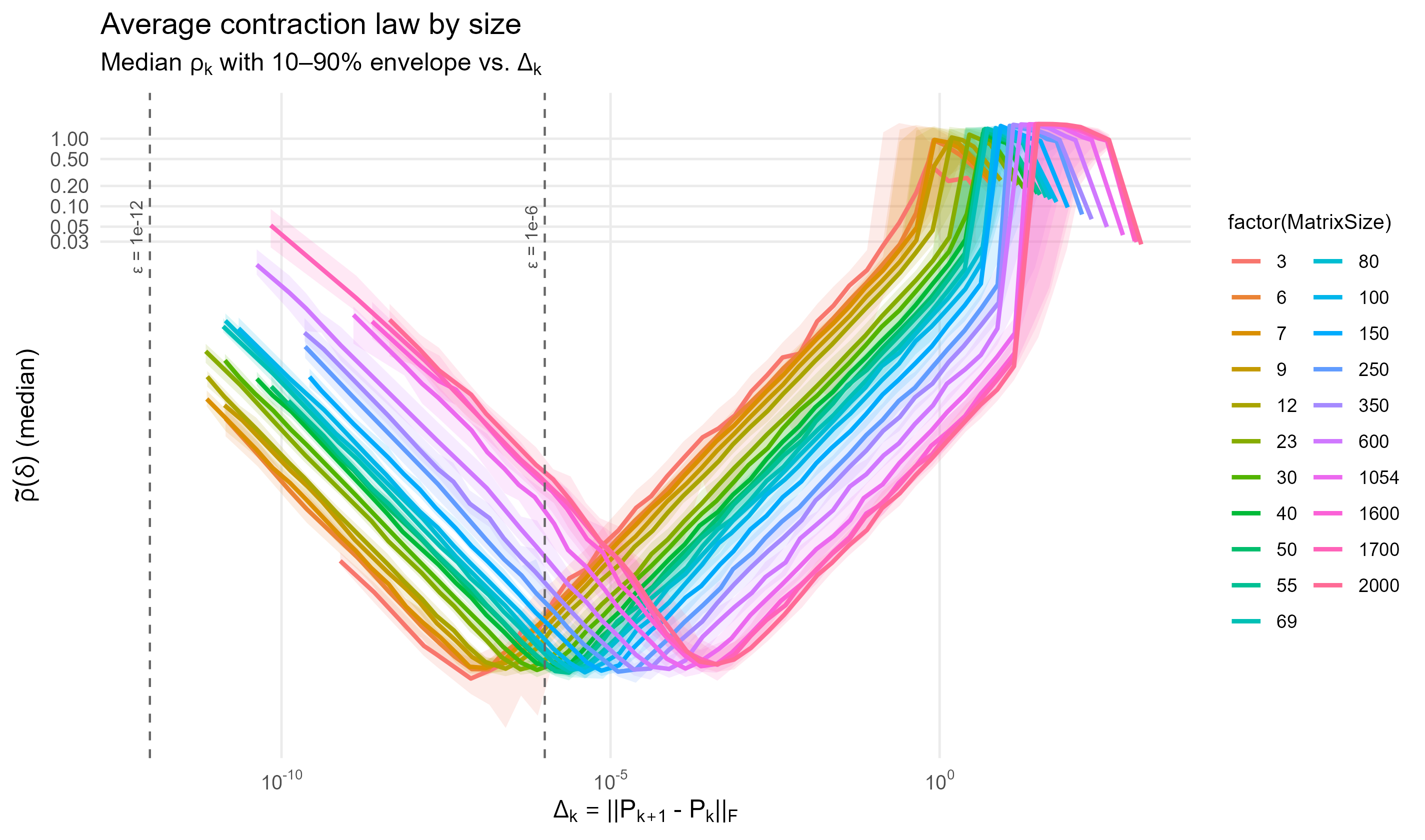}
  \caption{\textbf{Binned median contraction ratio by matrix size.}
  For each dimension \(n\), the curve shows the binned median
  \(\bar\rho_n(\delta_b)\) from~\eqref{eq:rho-bar}, and the shaded ribbon
  indicates the 10th–90th percentile range of \(\rho_k\) within each bin
  of \(\Delta_k\). For large \(\delta_b\), medians satisfy
  \(\bar\rho_n(\delta_b)\ll 1\) and the ribbon lies well below~\(1\),
  while for small \(\delta_b\) the medians approach \(1\) and the ribbon
  becomes a narrow band around~\(1\). The curves and ribbons for different
  \(n\) nearly coincide, showing that the observed \(\rho_k\)–\(\Delta_k\)
  behavior, including overshoots, is empirically dimension-independent.}
  \label{fig:ratio-by-size}
\end{figure}

\begin{figure}[!htbp]
  \centering
  \includegraphics[width=0.95\linewidth]{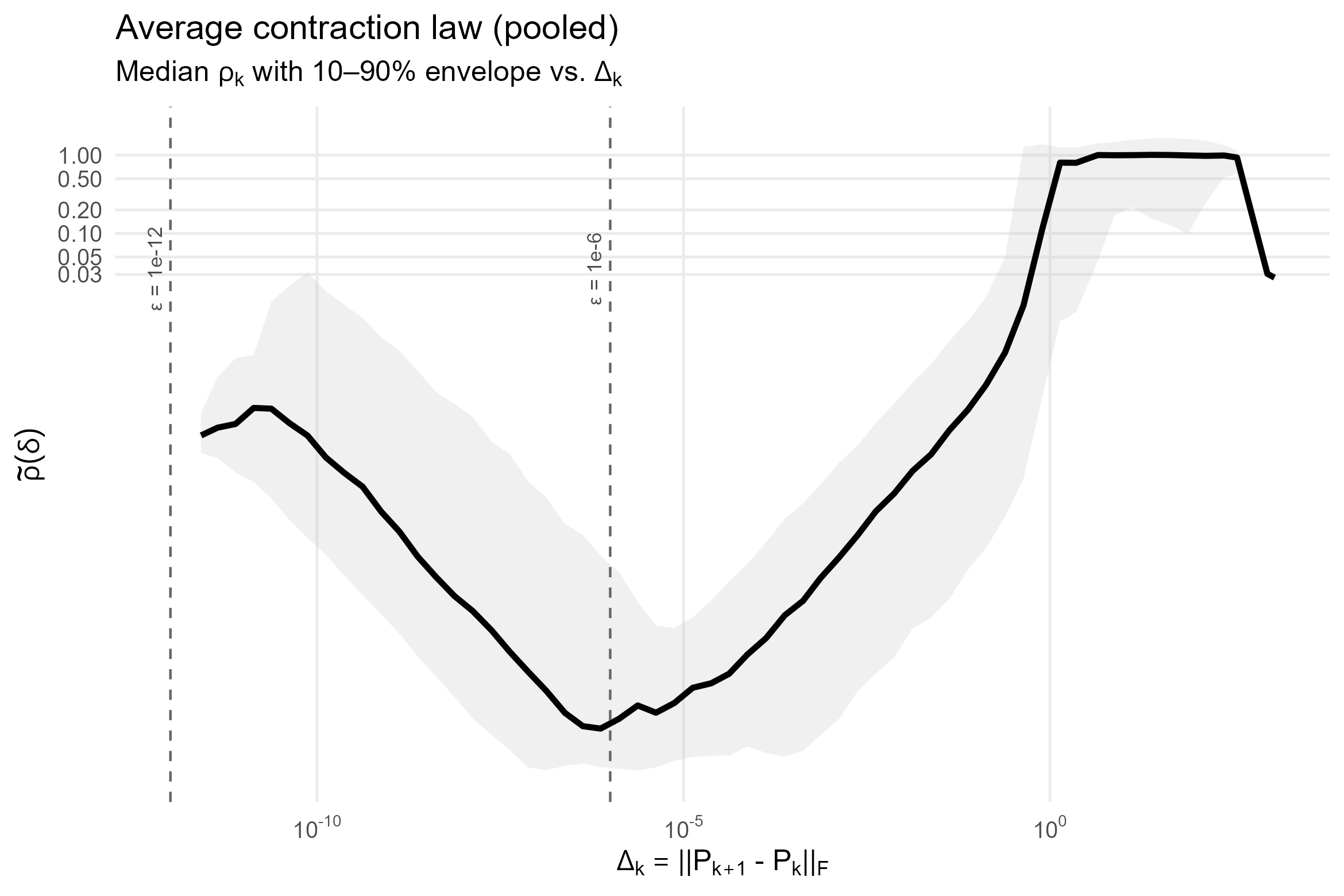}
  \caption{\textbf{Pooled binned median contraction law.}
  All pairs \((\Delta_k,\rho_k)\) from all sizes are pooled and binned in
  \(\Delta_k\). The curve shows the binned median contraction ratio
  \(\bar\rho(\delta_b)\), and the ribbon shows the 10th–90th percentile
  range. The pooled median and envelope reproduce the same V-shaped
  structure as in Figure~\ref{fig:ratio-by-size}: strong contraction at
  large \(\delta_b\), near-isometry with small overshoots at small
  \(\delta_b\), and a narrow spread around \(\rho=1\) in the small-step
  regime.}
  \label{fig:ratio-pooled}
\end{figure}

\paragraph{Interpretation and discussion.}
The organization of \((\Delta_k,\rho_k)\) into a V-shaped pattern, together
with the shrinking overshoot band as \(\Delta_k \downarrow 0\), is consistent
with the geometric structure described in Subsection~\ref{subsec:geometry},
but is established here only at the empirical level. Large step sizes
\(\Delta_k\) occur when the iterate is still far from satisfying the
row-wise mean-zero and unit-norm conditions; in this regime, centering and
normalization cause substantial changes to the rows, producing ratios
\(\rho_k \ll 1\). Once the rows lie close to the unit sphere in $\mathbf{1}^\perp$, the update moves mainly within the manifold $\Mn$, so $\Delta_{k+1}$ stays close to
$\Delta_k$ and $\rho_k \approx 1$, with only small overshoots above and below
unity. The resulting V-shape with a narrow overshoot band at small \(\Delta_k\) is
therefore a quantitative empirical signature of two regimes of the
iteration: a contractive regime when the rows lie far from the manifold $\Mn$, and a
near-isometric regime with small, bounded oscillations once they lie close to
it. This universal dimension-independent structure is one of the central
empirical findings of the study.

\subsection{Law IV: Uniformly bounded iteration counts}\label{subsec:law4}

For a fixed tolerance \(\varepsilon>0\), the convergence time is defined in
\eqref{eq:T-def} by
\[
T(n,\varepsilon)\;:=\;\min\Bigl\{k\in\mathbb{N}:\ \Delta_j<\varepsilon
\text{ for all } j\ge k\Bigr\},
\qquad
\Delta_k=\norm{P_{k+1}-P_k}_F.
\]

From \(N=1000\) independent trials at each size \(n\), samples
\(T_1(n,\varepsilon),\dots,T_N(n,\varepsilon)\) are obtained.

\begin{empiricallaw}[Bounded iteration counts]\label{law:law4}
For each fixed \(\varepsilon>0\), the distribution of \(T(n,\varepsilon)\)
remains in a tight, dimension-independent band. The median and upper
quantiles (including the 90\% and 95\% quantiles), the mean
\(\mu(n):=\mathbb{E}[T(n,\varepsilon)]\), and the standard deviation
\(\sigma(n):=\sqrt{\operatorname{Var}(T(n,\varepsilon))}\) are all uniformly
bounded in \(n\). Equivalently,
\[
\sup_n \operatorname{med}\bigl(T(n,\varepsilon)\bigr)<\infty,\quad
\sup_n Q_{0.95}\bigl(T(n,\varepsilon)\bigr)<\infty,\quad
\sup_n \mu(n)<\infty,\quad
\sup_n \sigma(n)<\infty.
\]
\end{empiricallaw}

For each dimension \(n\) and \(N=1000\) random initializations
\(P_0^{(t)}\), the first \(k\) satisfying \(\Delta_j^{(t)}<\varepsilon\) for
all \(j\ge k\) is recorded. The per-size statistics
\[
\begin{aligned}
\mu(n)&=\frac{1}{N}\sum_t T_t(n,\varepsilon),\;
\sigma(n)=\sqrt{\tfrac{1}{N-1}\sum_t (T_t(n,\varepsilon)-\mu(n))^{2}},\;
F_n(k)=\tfrac{1}{N}\sum_t\mathbf{1}\{T_t(n,\varepsilon)\le k\}.
\end{aligned}
\]

are then computed. The boundedness of \(\mu(n)\), \(\sigma(n)\), and high
quantiles of \(F_n(k)\) across all \(n\) supports \cref{law:law4}.

Empirically, Law~II shows that 
$\rho_k = \Delta_{k+1}/\Delta_k$ 
concentrates near $1$, with diminishing overshoot as $\Delta_k \to 0$. 
In practice this produces a small but persistent negative drift in 
$\log \Delta_k$, up to bounded fluctuations, which matches the uniformly 
bounded convergence times observed across all dimensions.

\begin{figure}
  \centering
  \includegraphics[width=0.95\linewidth]{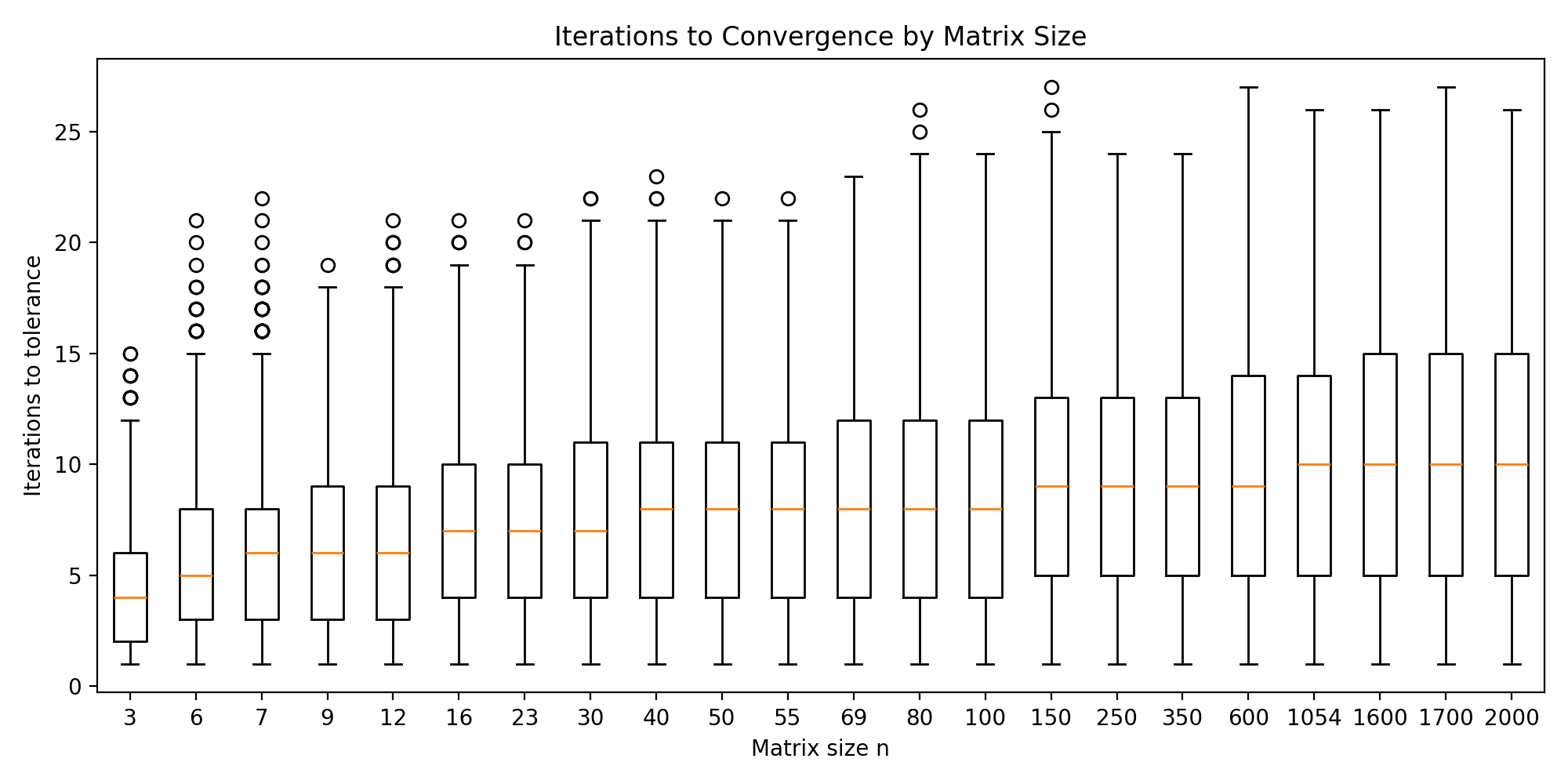}
  \caption{\textbf{Convergence time \(T(n,\varepsilon)\) across sizes (boxplots; \(N=1000\) trials per size).}
  For each \(n\), the box shows the interquartile range of
  \(T(n,\varepsilon)\), the central line is the median, whiskers mark the
  bulk, and points are outliers. The whisker band remains in a narrow range,
  here approximately between five and twenty iterations, with rare outliers
  below twenty-five. The width of this band does not grow with \(n\), which
  is direct distributional evidence that \(T(n,\varepsilon)\) is uniformly
  bounded in the matrix size.}
  \label{fig:iters-boxplot-law4}
\end{figure}

\begin{figure}[!htbp]
  \centering
  \includegraphics[width=0.95\linewidth]{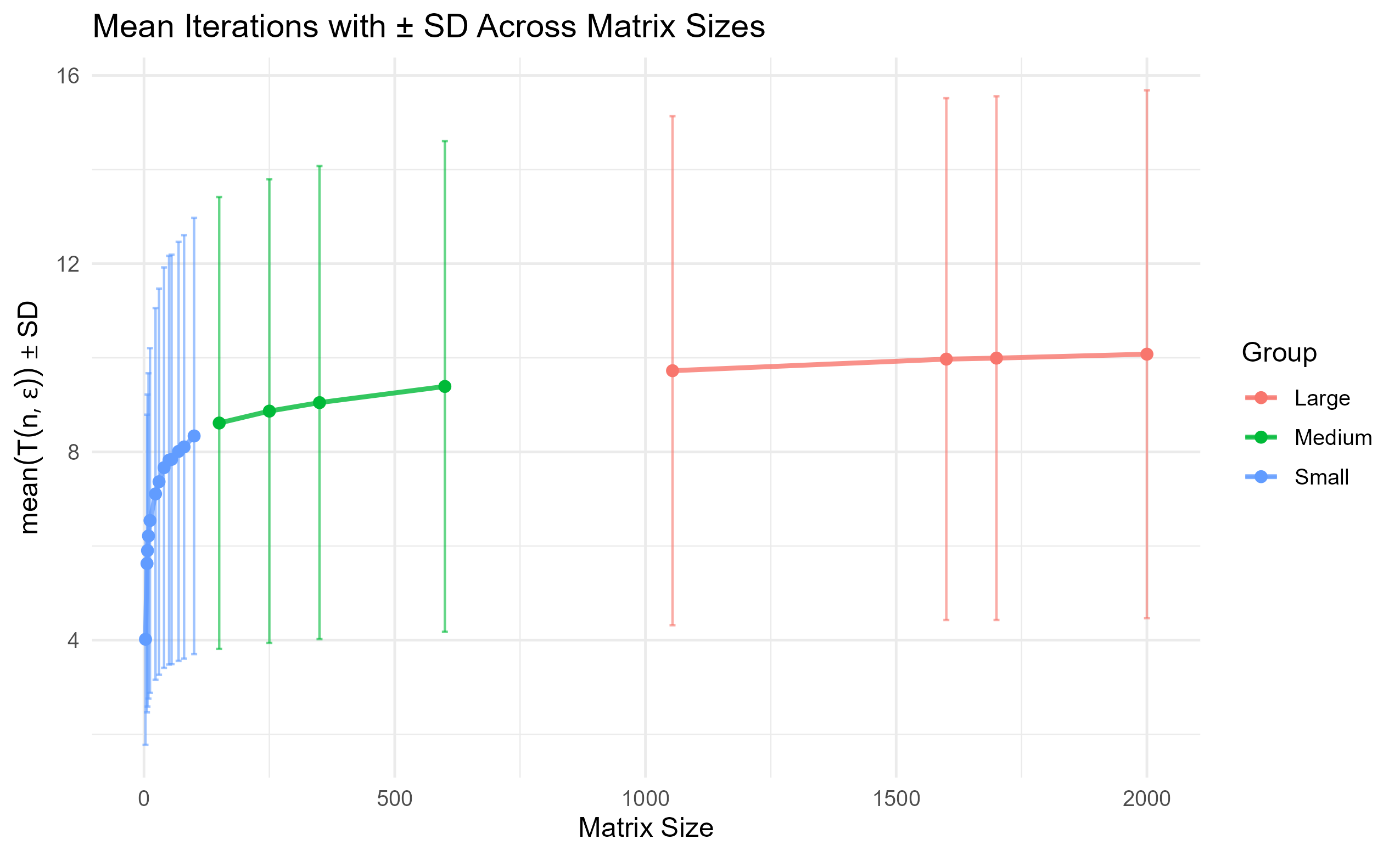}
  \caption{\textbf{Mean iterations \(\mu(n)\) with \(\pm\,\sigma(n)\) error bars, grouped by size.}
  From \(N=1000\) trials per size,
  \(
    \mu(n)=\tfrac{1}{N}\sum_{i=1}^N T_i(n,\varepsilon)
  \)
  and
  \(
    \sigma(n)=\sqrt{\tfrac{1}{N-1}\sum_{i=1}^{N}(T_i(n,\varepsilon)-\mu(n))^2}
  \).
  Both \(\mu(n)\) and \(\sigma(n)\) remain uniformly bounded and show no
  growth trend with \(n\). This second-moment view corroborates the
  boxplot evidence that \(T(n,\varepsilon)\) is of order one uniformly in
  \(n\).}
  \label{fig:mean-sd-law4}
\end{figure}

\begin{figure}[!htbp]
  \centering
  \includegraphics[width=0.95\linewidth]{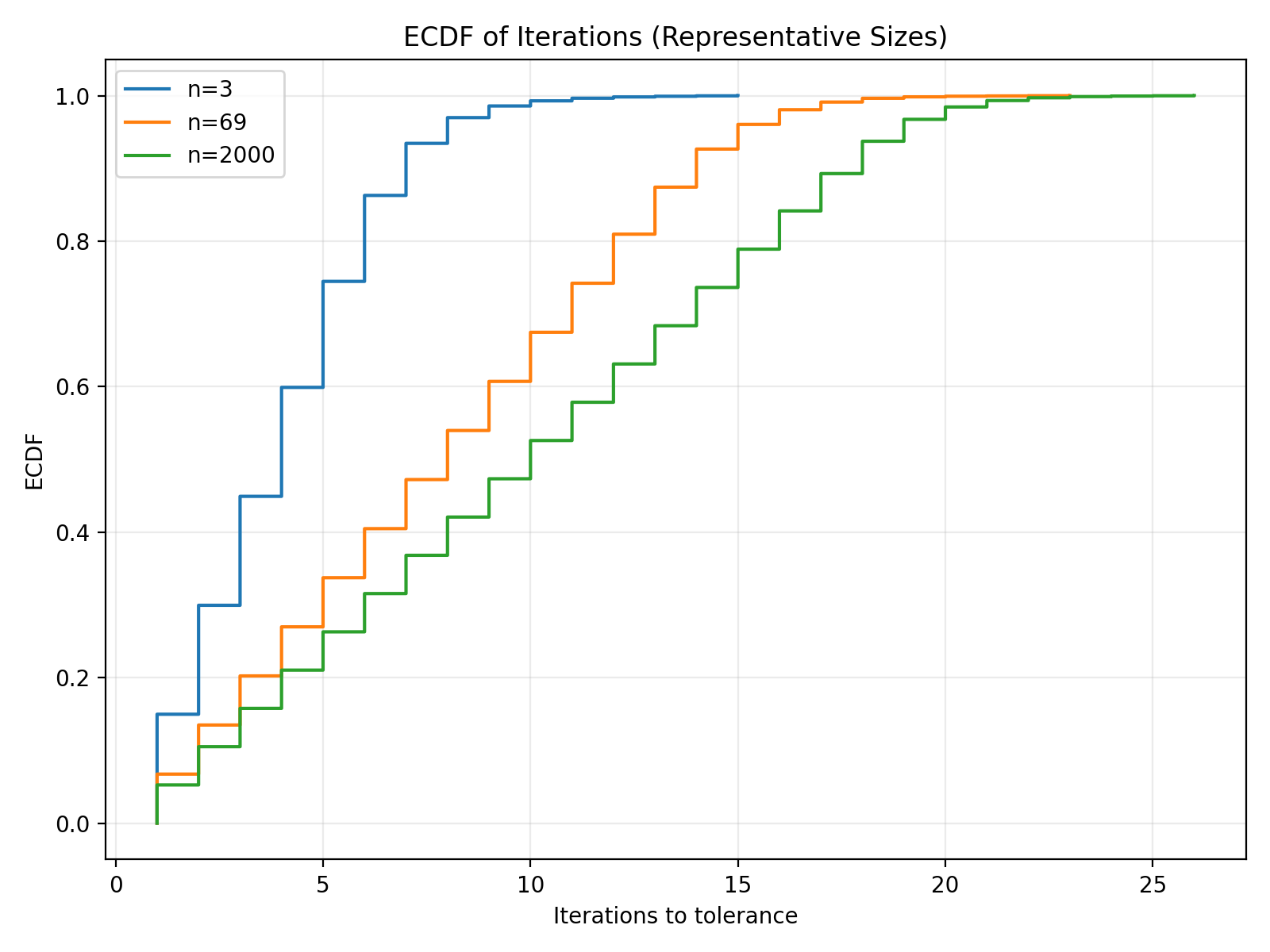}
  \caption{\textbf{Empirical distribution functions of iteration counts at representative sizes.}
  For each \(n\), the empirical distribution function
  \(F_n(k)=\tfrac{1}{N}\sum_{i=1}^{N}\mathbf{1}\{T_i(n,\varepsilon)\le k\}\)
  rises steeply and saturates rapidly. In the runs shown, the 90\% quantile
  occurs by \(k=7\) for \(n=3\), by \(k=15\) for \(n=69\), and by \(k=20\)
  for \(n=2000\). Almost all trajectories therefore terminate within a small
  number of steps, independent of dimension.}
  \label{fig:ecdf-reps-law4}
\end{figure}
\begin{table}[ht]
\centering
\small
\caption{Maximum observed iteration count $T_{\max}(n)$ over $N=1000$ trials per dimension 
for $\varepsilon = 10^{-12}$. The maximum number of steps never exceeds $28$ across the entire 
range $n\in[3,2000]$, demonstrating dimension-independent convergence time.}
\label{tab:iter-max}
\begin{tabular}{rr@{\qquad}rr@{\qquad}rr@{\qquad}rr}
\toprule
$n$ & $T_{\max}(n)$ &
$n$ & $T_{\max}(n)$ &
$n$ & $T_{\max}(n)$ &
$n$ & $T_{\max}(n)$ \\
\midrule
3    & 15 & 16   & 22 & 69   & 25 & 350  & 27 \\
6    & 21 & 23   & 23 & 80   & 25 & 600  & 27 \\
7    & 22 & 30   & 23 & 100  & 25 & 1054 & 28 \\
9    & 19 & 40   & 24 & 150  & 26 & 1600 & 28 \\
12   & 21 & 50   & 24 & 250  & 26 & 1700 & 28 \\
55   & 24 & 2000 & 28 &      &    &      &    \\
\bottomrule
\end{tabular}
\end{table}

As Table~\ref{tab:iter-max} shows, the iteration never exceeds 28 steps in any
of the 22{,}000 trials across dimensions $n\in[3,2000]$, providing definitive
numerical evidence that the convergence time is of order $O(1)$, entirely independent of dimension.

The uniform boundedness of the iteration count contrasts with many linear
scaling procedures such as Sinkhorn’s convex normalization
\cite{sinkhorn1964relationship}, whose iteration counts typically increase
with dimension.
The universality of \(T(n,\varepsilon)\) supports the finite boundedness of the cumulative variation defined in Equation \eqref{eq:total-variation}, whose value remains finite in all experiments. Moreover, while \(V_\infty\) grows approximately linearly with \(n\), its distribution at each fixed dimension is tightly concentrated around its median.


\section{Summary and Interpretation of the Empirical Laws}\label{sec:summary}

The four empirical laws established in \cref{sec:empirical} provide a
coherent description of the global behavior of the iteration
\(P_{k+1}=f(P_k)\). Each law is based on simple diagnostics yet together
they reveal a stable, dimension-independent structure across all tested
sizes.
This empirical stability was already observed in the classical clustering
literature \cite{mcquitty1968multiple,breiger1975clustering}, and implicitly
underlies Kruskal’s original CONCOR convergence report \cite{kruskalCONCOR}.

Law~I shows that the first update \(P_0\mapsto P_1\) produces a sharp
reduction in the Frobenius step size, with contraction ratio
\(\rho_0=\Delta_1/\Delta_0\) far below~\(1\) and essentially independent of
the matrix size. This reflects the fact that centering removes a large
common component of each row along \(\mathbf{1}\), and normalization enforces
unit length, jointly inducing a strong corrective motion toward the
mean-zero, unit-norm row structure imposed by the correlation operator.

Law~II shows that, after the initial drop, the sequence
\(\Delta_k\) decays in a nearly monotone fashion, with small overshoots that
remain uniformly bounded across dimensions and diminish as the iteration
progresses. Although no analytical
proof of finite total variation is known, the empirical behavior points
strongly in that direction. Indeed, across all dimensions, the total
variation
\[
V_\infty = \sum_{k\ge 0} \Delta_k
\]
is finite in every trial and for every trajectory and grows approximately linearly with \(n\), with narrow interquartile ranges at each dimension.

Law~III shows that the one-step contraction ratio \(\rho_k\) depends almost
exclusively on the current step size \(\Delta_k\). Large steps correspond to
strong contraction, while very small steps correspond to near-isometric
behavior. When pooled across all sizes and trials, the points
\((\Delta_k,\rho_k)\) align along the same characteristic V-shaped curve,
with a narrow band of overshoot around \(1\) in the small-step regime.
Once \(\Delta_k\) is specified, the dimension plays no further role.
The large-step branch reflects strong correction of deviations transverse to
the row-spherical manifold $\Mn$, whereas the small-step branch reflects the
nearly tangential motion of the iteration when the rows already lie close to
$\Mn$.

Law~IV shows that, for any fixed tolerance \(\varepsilon>0\), the number of
iterations \(T(n,\varepsilon)\) required for stabilization remains tightly
bounded across all matrix sizes. This is consistent with the eventual
negative drift of \(\log \Delta_k\) suggested by Law~II and with the
universal \(\rho\)--\(\Delta\) relationship of Law~III. Empirically, almost
all trajectories stabilize in a small, dimension-independent number of
steps, even for the largest matrices tested.

Taken together, the four laws reveal a geometric picture of the composite
map \(f=\Phi\circ\Psi\). The centering–normalization step pushes the iterate
sharply toward the underlying row-spherical geometry and heavily damps
deviations transverse to it. Once the rows lie close to this geometry, the
dynamics become almost tangential, producing small, stable steps whose
behavior is accurately captured by the universal V-shaped contraction law.
The Gram map transfers these changes to the correlation matrix without
introducing dimension-dependent distortions, explaining the strong initial
drop, the slow and controlled decay that follows, and the
dimension-independent behavior of all quantitative measures.

A closer examination of the cumulative step length
$V_\infty = \sum_{k\ge0}\Delta_k$ reveals an additional geometric feature of
the iteration.  
For every tested dimension, $V_\infty$ is finite, in agreement with the
near-monotone decay of $\Delta_k$ and the bounded overshoot behavior of
Law~II.  
However, unlike the iteration count \(T(n,\varepsilon)\) of Law~IV, the total
variation is \emph{not} dimension-independent.  
Empirically, \(V_\infty\) grows with \(n\) in an approximately linear fashion, \(V_\infty \approx C n\) with \(C \approx 2.3\), while remaining tightly concentrated around its median at each dimension.

This behavior reflects the geometry of the row-spherical manifold
\[
\Mn=\{Z\in\mathbb{R}^{n\times n}:\ \mathbf{1}^\top Z(i,:)=0,\ 
\|Z(i,:)\|_2=1\ \text{for all } i\},
\]
onto which the operator $\Psi$ maps each centered row of $P_k$.  
For large $n$, the initial matrix $P_0$ typically lies farther from $\Mn$ in
Frobenius distance, because its rows deviate more substantially—in both mean
and magnitude—from the unit sphere contained in $\mathbf{1}^\perp$.  
The early iterations therefore traverse a longer path in the ambient space.
Once the iterate approaches $\Mn$, however, the dynamics enter the
near-isometric regime documented in Law~III, and the number of subsequent
corrective steps remains uniformly bounded in $n$.

Thus, the scaling of \(V_\infty\) with dimension reflects a clear geometric
separation: the \emph{distance} from the initial state to the row-spherical
manifold \(\Mn\) increases with \(n\), but the \emph{number of steps}
required to reach and remain near \(\Mn\) does not.
This reconciles the growth of \(V_\infty\) with the striking
dimension-independence seen in Laws~III and~IV.

The empirical laws also clarify why convergence has been consistently
observed since the early work of McQuitty, Breiger--Boorman--Arabie, and
Chen. The map is strongly contracting in directions that alter row means or
row norms, and nearly isometric along directions tangent to the mean-zero,
unit-norm row manifold. This dual behavior is straightforward to verify locally but challenging to analyze globally due to nonlinearity and the structure of the fixed-point set of block-$\{\pm1\}$ patterns. Despite these difficulties, the empirical laws
demonstrate that the iteration behaves regularly across all dimensions and
admits stable patterns of contraction, overshoot control, and convergence
time.

The main theoretical challenge is to convert these empirical regularities
into uniform estimates valid for all initial matrices. Law~II suggests a
Lyapunov-type structure with finite total variation, and Law~III provides a
precise functional dependence between \(\Delta_k\) and \(\rho_k\) that any
future theory must explain. Establishing such properties analytically would
close a long-standing gap in the understanding of iterated correlation maps
and complete the convergence picture initiated in the earlier literature.

\subsection{Why global convergence has remained an open problem}\label{subsec:why-open}

The empirical evidence presented in this work—four sharp, dimension-independent regularities that hold uniformly up to dimension $n=2000$—reveals a degree of order that appears almost too good to be true for a nonlinear iteration. Yet, despite more than fifty years of consistent empirical success across psychometrics, social-network analysis, and data visualization, no proof of global convergence from arbitrary initial conditions has ever been published. At first sight, standard fixed-point theorems might seem applicable.
Brouwer’s theorem guarantees that the continuous map $f$ has at least one
fixed point on the compact, convex elliptope $\mathcal{E}_n$. Moreover, if a global
contraction bound
\[
\|f(P)-f(Q)\|_F \le c\,\|P-Q\|_F \qquad 
\]
held uniformly for all $P,Q$ with some \(c<1\), then Banach’s contraction principle would
immediately imply uniqueness of the fixed point and global convergence of the
iteration $P_{k+1}=f(P_k)$. However, none of the hypotheses needed for such
arguments are satisfied. The structural obstacles are decisive:

\begin{enumerate}
\item \emph{The fixed-point set is finite, and determined by discrete block structure.}  
Kruskal (ca.\ 1965, as summarized in~\cite{breiger1975clustering}) established that the only fixed points reachable by the iteration are block-$\{\pm1\}$ correlation matrices: symmetric positive semidefinite matrices with unit diagonal whose off-diagonal entries are exactly $+1$ within blocks and $-1$ between blocks (up to row/column permutation and sign flips of entire blocks).  
This set $\mathcal{F}_n \subset \mathcal{E}_n$ is finite for each $n$, within the class of iterates reachable from nondegenerate initial conditions, and it grows exponentially with $n$. Although it lies inside the convex elliptope $\mathcal{E}_n$ of all correlation matrices, $\mathcal{F}_n$ itself, was not fully characterized by \cite{chen2002gap,kruskalCONCOR}. Consequently, virtually all global convergence frameworks that rely on convexity of the target set—Fejér monotonicity, proximal methods, or Bregman projections~\cite{bauschke2017convex}—are inapplicable.
A more systematic analysis of these algebraic fixed points is developed in the
forthcoming work ~\cite{zusmanovich2025fixedpoints}, but a complete
characterization for general \(n\) remains open.

\item \emph{The map is nonlinear and not uniformly contractive.}  
Law~III demonstrates a near-isometric regime for small perturbations, so Banach’s contraction principle cannot be applied globally.
\item \emph{The iteration is not a standard Mann-type averaging scheme.}
Classical mean-value iteration methods \cite{mann1953mean} study convex
combinations of a point with its image under a nonexpansive map. The
correlation iteration applies the full nonlinear map $f$ at each step
without averaging and, moreover, lacks global nonexpansivity, so
Mann-type convergence results cannot be invoked directly.

\item \emph{Contraction is strongly anisotropic.}  
The iteration contracts aggressively in directions that alter row means or
row norms (Law~I and the transverse branch of Law~III), but is nearly
isometric along directions tangent to the mean-zero, unit-norm row
manifold (the tangential branch of Law~III). No known metric turns the map
into a uniform contraction on the whole space.

\item \emph{Sinkhorn-type arguments fail.}  
Sinkhorn–Knopp scaling converges globally because the target set (doubly-stochastic matrices) is convex and the map is a contraction in Hilbert’s projective metric~\cite{sinkhorn1964relationship}. Here the effective target is the discrete set $\mathcal{F}_n$, and no analogous metric exists.
\end{enumerate}

Kruskal’s local linearization shows that, once sufficiently close to any individual block-$\{\pm1\}$ fixed point, convergence to that exact point is geometrically fast with spectral radius strictly smaller than $1$~\cite{breiger1975clustering}. However, the global trajectory from a generic initial matrix must cross regions far from all fixed points, and the basins of attraction of the exponentially many fixed points remain unknown.

The four empirical laws reported here are therefore all the more remarkable: they reveal that the dynamics are extraordinarily regular and essentially dimension-independent \emph{despite} the severe structural barriers above. They constitute the first quantitative, global portrait of the correlation iteration and provide the precise, testable benchmarks that any future analytical convergence proof must reproduce.

\section{Discussion and open problems}\label{sec:openproblems}

The empirical laws presented above clarify several structural properties of
the iteration but also highlight what remains to be proved. Several open
mathematical problems are listed here.

\paragraph{Global convergence.}
Classical local arguments support near-pattern geometric convergence, but a
complete convergence theorem valid for all initializations and all
dimensions is not available.

\paragraph{Finite total variation.}
Empirically, the sum of step sizes
\(\sum_k \Delta_k\) is finite for every trajectory and grows approximately linearly with the dimension, \(V_\infty \approx C n\). A proof would require
showing that overshoots remain uniformly bounded and that the effective
contraction ratio eventually drops below unity on average, together with a global upper bound of order \(O(n)\).

\paragraph{Uniform control of overshoots.}
Law~II shows that overshoots are small and bounded across all sizes. A
theoretical explanation of this phenomenon is missing and would likely
require a detailed study of the geometry of the centering and normalization
maps.

\paragraph{Universal contraction law.}
Law~III identifies a dimension-independent relation between the step size
\(\Delta_k\) and the ratio \(\rho_k\). Understanding why this V-shaped
profile emerges and why it is independent of \(n\) is an open analytical
question.

\paragraph{Characterization of fixed points.}
The fixed--point set consists of block-\(\{\pm 1\}\) correlation patterns.
A detailed analysis of these fixed points, including stability properties and basins of attraction, is developed in the preprint~\cite{zusmanovich2025fixedpoints}.

\paragraph{Lyapunov functions.}
The iteration appears to admit a Lyapunov structure with strong contraction
in directions orthogonal to the row-spherical manifold $\mathcal{M}_n$ and
near-isometry along directions tangent to $\mathcal{M}_n$, together with
finite cumulative change. Constructing such a function explicitly remains an
open challenge.

\paragraph{Dimension-independence.}
Many observed properties---first-step contraction, overshoot bounds,
convergence time—are stable in \(n\). Explaining this independence
analytically is a key theoretical question.

The empirical laws identified here therefore supply the structural benchmarks that any full analytical theory must reproduce, and they delineate a precise pathway toward resolving a problem that has remained open since the 1960s.


\begin{thebibliography}{99}

\bibitem{alhajjhassan2025framework}
I.~AlHajj Hassan.
\newblock \emph{Computational Framework for Experiments on Iterated Correlation Matrices}.
\newblock Zenodo, 2025.
\newblock \href{https://doi.org/10.5281/zenodo.17794063}{https://doi.org/10.5281/zenodo.17794063}.

\bibitem{zusmanovich2025fixedpoints}
I.~AlHajj Hassan and P.~Zusmanovich.
\newblock ``Fixed points of iterated correlation matrices.''
\newblock Preprint, 2025.

\bibitem{bauschke2017convex}
H.~H. Bauschke and P.~L. Combettes.
\newblock \emph{Convex Analysis and Monotone Operator Theory in Hilbert Spaces},
2nd edition.
\newblock Springer, 2017.
\newblock \href{https://doi.org/10.1007/978-3-319-48311-5}{https://doi.org/10.1007/978-3-319-48311-5}.

\bibitem{breiger1975clustering}
R.~L. Breiger, S.~A. Boorman, and P.~Arabie.
\newblock An algorithm for clustering relational data with applications to social
network analysis and comparison with multidimensional scaling.
\newblock \emph{Journal of Mathematical Psychology}, 12(3):328--383, 1975.
\newblock \url{https://doi.org/10.1016/0022-2496(75)90028-0}{https://doi.org/10.1016/0022-2496(75)90028-0}.

\bibitem{chen2002gap}
C.-C. Chen.
\newblock Generalized association plots: Information visualization via iteratively
generated correlation matrices.
\newblock \emph{Statistica Sinica}, 12(1):7--29, 2002.
\newblock \url{http://www3.stat.sinica.edu.tw/statistica/}{http://www3.stat.sinica.edu.tw/statistica/}.

\bibitem{goebel1990topics}
K.~Goebel and W.~A. Kirk.
\newblock \emph{Topics in Metric Fixed Point Theory}.
\newblock Cambridge University Press, 1990.
\newblock \url{https://doi.org/10.1017/CBO9780511526152}{https://doi.org/10.1017/CBO9780511526152}.

\bibitem{huang2019iterative}
L.~Huang, D.~Yang, and B.~Lang.
\newblock Iterative normalization: Beyond standardization towards efficient whitening.
\newblock In \emph{Proceedings of the IEEE Conference on Computer Vision and
Pattern Recognition (CVPR)}, pages 4874--4883, 2019.
\newblock \url{https://openaccess.thecvf.com/content_CVPR_2019/html/Huang_Iterative_Normalization_Beyond_Standardization_Towards_Efficient_Whitening_CVPR_2019_paper.html}{OpenAccess: CVPR 2019}.

\bibitem{kruskalCONCOR}
J.~B. Kruskal.
\newblock A theorem about CONCOR.
\newblock Technical Report MH~2C--571, Bell Laboratories, Murray Hill, NJ, 1978.

\bibitem{laurent1996elliptope}
M.~Laurent and S.~Poljak.
\newblock On the facial structure of the set of correlation matrices.
\newblock \emph{SIAM Journal on Matrix Analysis and Applications},
17(2):530--547, 1996.
\newblock \url{https://doi.org/10.1137/S089547989325030X}{https://doi.org/10.1137/S089547989325030X}.

\bibitem{mann1953mean}
W.~R. Mann.
\newblock Mean value methods in iteration.
\newblock \emph{Proceedings of the American Mathematical Society},
4(3):506--510, 1953.
\newblock \url{https://doi.org/10.1090/S0002-9939-1953-0054846-3}{https://doi.org/10.1090/S0002-9939-1953-0054846-3}.

\bibitem{mcquitty1968multiple}
L.~L. McQuitty.
\newblock Multiple clustering revisited: Comments, comparisons, new approaches.
\newblock \emph{Multivariate Behavioral Research}, 3(4):431--479, 1968.
\newblock \url{https://doi.org/10.1207/s15327906mbr0304_1}.


\bibitem{opial1967weak}
Z.~ Opial.
\newblock Weak convergence of the sequence of successive approximations for
nonexpansive mappings.
\newblock \emph{Bulletin of the American Mathematical Society}, 73(4):591--597, 1967.
\newblock \url{https://doi.org/10.1090/S0002-9904-1967-11761-0}{https://doi.org/10.1090/S0002-9904-1967-11761-0}.

\bibitem{saunderson2012elliptope}
J.~Saunderson, V.~Chandrasekaran, P.~A. Parrilo, and A.~S. Willsky.
\newblock Diagonal and low-rank matrix decompositions, correlation matrices, and
ellipsoid fitting.
\newblock \emph{SIAM Journal on Matrix Analysis and Applications},
33(3):813--839, 2012.
\newblock \url{https://doi.org/10.1137/110837851}{https://doi.org/10.1137/110837851}.

\bibitem{sinkhorn1964relationship}
R.~Sinkhorn.
\newblock A relationship between arbitrary positive matrices and doubly stochastic
matrices.
\newblock \emph{The Annals of Mathematical Statistics},
35(2):876--879, 1964.
\newblock \url{https://doi.org/10.1214/aoms/1177703599}{https://doi.org/10.1214/aoms/1177703599}.

\end{thebibliography}
\end{document}